\documentclass{article}
\usepackage{amsmath,amssymb}
\begin{document}


\newcommand{\Kfrac}[2]{\ensuremath{\displaystyle\frac{{#1}}{{#2}}}}
\newcommand{\cpl}{\ensuremath{\mathbb{C}}}
\newcommand{\eul}{\ensuremath{\mathrm{i}}}
\newcommand{\reel}{\ensuremath{\mathbb{R}}}
\newcommand{\dif}[1]{\ensuremath{\:\mathrm{d}{#1}}}
\newcommand{\dpar}[1]{\ensuremath{\partial_{{#1}}}}
\newcommand{\nm}[1]{\ensuremath{\|{#1}\|}}
\newcommand{\ps}[2]{\ensuremath{<\!\!\!~{#1}\,|\,{#2}~\!\!\!>}}
\newcommand{\ex}[1]{\ensuremath{\mathrm{e}^{{#1}}}}
\newcommand{\exni}[1]{\ensuremath{\ex{-\eul {#1}}}}
\newcommand{\ph}{\ensuremath{\varphi}}
\newcommand{\kfrac}[2]{\ensuremath{{{#1}}/{{#2}}}}
\newcommand{\exi}[1]{\ensuremath{\ex{\eul {#1}}}}
\newcommand{\Inv}[1]{\ensuremath{\frac{1}{{#1}}}}
\newcommand{\ens}[2]{\ensuremath{\{{#1}\;;\;{#2}\}}}
\newcommand{\Ens}[2]{\ensuremath{\left\{{#1}\;;\;{#2}\right\}}}
\newcommand{\KInv}[1]{\ensuremath{\displaystyle\Inv{{#1}}}}
\newcommand{\eps}{\ensuremath{\varepsilon}}
\newcommand{\sgn}[1]{\ensuremath{\mathrm{sgn}({#1})}}
\newcommand{\Klim}{\ensuremath\displaystyle\lim}
\newcommand{\rel}{\ensuremath{\mathbb{Z}}}
\newcommand{\sgl}[1]{\ensuremath{\{{#1}\} }}
\newcommand{\TroisEtoiles}{\mbox{}\\ \begin{center} * * * \end{center} \mbox{}\\ }
\newcommand{\kTroisEtoiles}{\begin{center} * * * \end{center}}
\newcommand{\I}{\ensuremath{\mathit{I}}}
\newcommand{\II}{\ensuremath{\mathit{II}}}
\newcommand{\III}{\ensuremath{\mathit{III}}}
\newcommand{\dist}[2]{\ensuremath{\mathrm{dist}({#1},{#2})}}
\newcommand{\kInv}[1]{\ensuremath{\displaystyle 1/{#1}}}
\newcommand{\Kprod}{\ensuremath{\displaystyle\prod}}
\newcommand{\Ksum}{\ensuremath{\displaystyle\sum}}
\newcommand{\dem}{\ensuremath{\frac{1}{2}}}
\newcommand{\adh}[1]{\ensuremath{\overline{{#1}}}}
\newcommand{\vp}{\ensuremath{\mathrm{p.v.}\int}}
\newcommand{\Kint}{\ensuremath{\displaystyle\int}}
\newcommand{\nat}{\ensuremath{\mathbb{N}}}
\newcommand{\fml}[3]{\ensuremath{\lbrace {#1} \rbrace_{{#2}}^{{#3}}}}
\newcommand{\dbl}[2]{\ensuremath{\stackrel{\scriptstyle {#1}}{{#2}}}}

\newtheorem{thm}{Theorem}

\newenvironment{exe}{\mbox{}\\ \noindent \textbf{Example}\ }{\mbox{}\\}

\newenvironment{kcor}{\mbox{}\\ \noindent \textbf{Corollary}\itshape\ }

\newenvironment{pr}{\noindent \textbf{Proof}\ }{\mbox{}\hfill $\Box$\\}

\newtheorem{prop}{Proposition}

\newenvironment{rem}{\noindent \textbf{Remark}\itshape\ }{\mbox{}\\}

\newenvironment{cor}{\mbox{}\\ \noindent \textbf{Corollary}\itshape\ }{\mbox{}\\}

\newenvironment{sch}{\noindent \textbf{Scholium}\itshape\ }{\mbox{}\\}

\newenvironment{krem}{\noindent \textbf{Remark}\itshape\ }

\newenvironment{lem}{\mbox{}\\ \noindent \textbf{Lemma}\itshape\ }{\mbox{}\\}


\newcommand{\notsubseteq}[1]{\subseteq\hspace{-#1}/\ }

 \newcommand{\comm}[1]{\textsc{[{#1}]}}
 \renewcommand{\comm}[1]{}

\newcommand{\Rtr}{\mbox{}\\}

\newenvironment{notes}{\Rtr\em\noindent{Notes:}\begin{enumerate}}{\end{enumerate}}

 \newcommand{\Hil}{\mathcal{H}}
 \newcommand{\conj}[1]{\overline{{#1}}}
 \renewcommand{\ps}[2]{\langle{#1},{#2}\rangle}
 \newcommand{\cutoff}{\chi}
 \newcommand{\str}{\mathcal{S}}
 \newcommand{\cjh}[1]{\widetilde{{#1}\,}}
 \renewcommand{\mod}{\mathrm{mod}\,}
 \newcommand{\pw}{\mbox{\textsl{pw}}}
 \newcommand{\reg}{\mathcal{R}}
 \newcommand{\Alpha}{\mathcal{A}}
 \newcommand{\Nev}{\mathcal{N}}

\title{Weighted Paley-Wiener spaces and mountain chain axioms: a detailed exposition}
\author{Philippe Poulin}
\maketitle
\abstract{%
The present work reviews the second half of Lyubarskii and Seip's paper, \emph{Weighted Paley--Wiener Spaces}.
Axioms defining a larger class of de Branges spaces are abstracted, allowing us to state and prove their results at a higher level of generality. 
} 

\section{Introduction}

In 2002, in the exploratory paper \cite{LyuSeip} Lyubarskii and Seip have introduced the notion of \emph{weighted Paley--Wiener space}.
Such a space consists of entire functions whose norms are comparable with the $L^2$-norms against $M(x)^{-2}\dif{x}$, where $M(x)$, the so-called 
 \emph{majorant}, is the norm of the reproducing kernel at $x\in\reel$ (see Section~\ref{Sec:WPWS} for complete definitions).

By a profound analysis, the authors obtained a concrete characterization of all weighted Paley--Wiener spaces. 
Each of them may be viewed as a perturbation of a classical Paley--Wiener space $L^2_a$, obtained by moving slightly the zeroes of its
 associated Hermite--Biehler function $\sin(a(z+\eul))$, so the distance of two successive zeroes remains comparable with $1$.
The explicit form of these spaces is similar to $L^2_a$, but with an ``exponential type condition" now involving a certain potential on the real line.

In the first part of their analysis the authors studied the distribution of the zeroes of an arbitrary Hermite--Biehler function associated with
 a weighted Paley--Wiener space.
They deduced that the derivative of its phase must satisfy properties making it comparable with a ``mountain chain". 
They then studied all de Branges spaces with these mountain chain properties, and specialized their results to weighted
 Paley--Wiener spaces.

The present expository paper shall review the second part of Lyubarskii and Seip's analysis, namely, how to obtain from the mountain chain properties:
\begin{itemize}
 \item
  The explicit form of the majorant of an arbitrary weighted Paley-Wiener space;
 \item
  The concrete form of an arbitrary weighted Paley--Wiener space;
 \item
  The properties to add to the mountain chain axioms so the resulting space is a given weighted Paley--Wiener space.
\end{itemize}

As far as possible, Lyubarskii and Seip's work is stated at its highest level of generality, in terms of spaces satisfying the mountain
 chain properties. 
These last properties are listed in Section \ref{Sec:MCA}, but the proof that all weighted Paley--Wiener spaces must satisfy them is not reviewed:
 readers are thus invited to refer to the lemma~1 and its proof in \cite{LyuSeip}.

Theorem~2 in \cite{LyuSeip} represented more difficulty: it provides a concrete realization of weighted Paley--Wiener spaces, but it also tells more,
 by characterizing all weighted Paley--Wiener spaces sharing the same majorant on the real line. 
The proof of this last part could not be completed without introducing new ideas (about the uniform density), and hence has been the object of a
 research paper apart \cite{Poulin}; it shall not be commented here.
Considerations about sampling and interpolation are also omitted, except the ones used for proving Theorem~4 in the original paper.

\paragraph{Paper's outline}
Basic definitions in de Branges theory and the formal definition of weighted Paley--Wiener space are given in Section~\ref{Sec:WPWS}.
Concrete examples of weighted Paley--Wiener spaces and their link with potential theory are then presented in Section~\ref{Sec:APE}.
Notice that these last examples cover indeed all possible cases, as proven in a subsequent section (Theorem~\ref{Thm:GFWPWS:GFWPWS}).
The mountain chain properties used in Lyubarskii and Seip's analysis are then stated as axioms in Section~\ref{Sec:MCA}, yielding the more
 general notion of \emph{MC-space}.
Then, the three aforementioned main results are proven in their respective sections.

The present paper aims to present Lyubarskii and Seip's analysis in more detail and state their results in their full generality.
For the sake of clarity steps are presented in a different order and new terminology is introduced. 
Needless to say, huge efforts were done for avoiding imprecision and correcting the unavoidable typos.
It is hoped that, for such a technical work, the present complement to \cite{LyuSeip} will be appreciated.

\paragraph{Acknowledgements}

\paragraph{Notation} We shall denote the derivative of $f$ with respect to $x$ by $\dpar{x}f$. 
Also, as in the original paper we shall use the following notation: $f\lesssim g$ denotes the existence of a constant $C$ such that $f\leq Cg$,
 while $f\simeq g$ means that $f$ and $g$ are comparable ($f\lesssim g$ and $g\lesssim f$).

\section{Weighted Paley--Wiener spaces}\label{Sec:WPWS}

A Hilbert space $\Hil$ of entire functions is a \emph{de Branges space} if it satisfies the following axioms \cite[th.23]{deBranges}:
\emph{%
 \begin{enumerate}
  \item
   The linear functional $\Hil\to\cpl,\ f\mapsto f(z_0)$ is bounded for all $z_0\in\cpl$;
  \item
   If $f(z)\in\Hil$, then $f^*(z)$ also belongs to $\Hil$ and has the same norm as $f(z)$;
  \item
   If $f(z)\in\Hil$ and $f(z_0)=0$, then $f(z)\Kfrac{z-\conj{z_0}}{z-z_0}$ also belongs to $\Hil$ and has the same norm as $f(z)$.
 \end{enumerate}
}
\noindent By the first axiom, $\Hil$ admits a \emph{reproducing kernel}, that is, a function $k_w(z)$ of the variables $w,z\in\cpl$ such that: $k_w\in\Hil$ for all $w\in\cpl$, and
  $$\ps{f}{k_w}_{\Hil}=f(w)\mbox{ \ for all }f\in\Hil.$$
The \emph{majorant} of $\Hil$ at $z\in\cpl$ is then defined as
 $$M(z)=\nm{k_z}_\Hil=\sup_{\nm{f}_\Hil=1}|f(z)|.$$

Let $M(x)$ be the restriction of $M$ to the real axis.
Following Lyubarskii and Seip \cite{LyuSeip}, we shall say that $M(x)$ is a \emph{majorant-weight} if
 \begin{enumerate}
  \item
   $M(x)>0$ for all $x\in\reel$;
  \item
   $\nm{f}_{\Hil}\simeq\nm{f/M}_2$ for all $f\in\Hil$.
 \end{enumerate}
Then, the corresponding $\Hil$ is called a \emph{weighted Paley--Wiener space}.\\

\begin{rem}
To be a weighted PW-space is invariant under equality with norm equivalence.
However, if an isometry maps a weighted PW-space to a de Branges space, this last is not necessarily
 a weighted PW-space.
\end{rem}

Let us recall useful facts about de Branges spaces.
An \emph{Hermite--Biehler function} $E$ is an entire function satisfying $|E(z)|>|E(\conj{z})|$ for all $z\in\cpl^+$.
Such a function may be factorized as \cite[th.V.6]{LevinOld}
 \begin{equation}\label{Eqn:WPWS:FactoHermiteBiehler}
  E(z)=C z^m \ex{h(z)}\exni{\alpha z}\prod_{\lambda\in\Lambda}(1-z/\lambda)\ex{p_\lambda(z)},
 \end{equation}
where $C\in\cpl$, $h(z)$ is real-entire, $\alpha\geq 0$, $\Lambda$ is a family of nonzero elements lying in the closed
 lower half-plane (with possible repetitions) and for all $\lambda$, $p_\lambda(z)$ is a polynomial with real coefficients.
Conversely, given such a factorization (where $\alpha\neq 0$ or $\Lambda\notsubseteq{3.5mm}\reel$), if the
 right-hand side in~(\ref{Eqn:WPWS:FactoHermiteBiehler}) defines an entire function, then it is in the
 Hermite--Biehler class.%
\footnote{%
Indeed, for $\Im z>0$, $|\ex{h(\bar{z})}|=|\ex{h(z)}|$ since $h$ is real-entire,
 $|\exni{\alpha\bar{z}}|\leq|\exni{\alpha z}|$ since $\alpha\geq 0$, $|1-\bar{z}/\lambda|\leq|1-z/\lambda|$ since $\Im\lambda\leq 0$,
 and $|\ex{p_\lambda(\bar{z})}|=|\ex{p_\lambda(z)}|$ since the coefficients of $p_\lambda$ are real.}

In the case where $E$ does not have real zeroes, its restriction to the real axis may be written \cite[prob.48]{deBranges}
 $$E(x)=|E(x)|\exni{\ph(x)},$$
where the \emph{phase}, $\ph(x)$, is real-analytic and well-defined (up to the addition of $2k\pi$).
The factorization~(\ref{Eqn:WPWS:FactoHermiteBiehler}) then implies
 \begin{equation} \label{Eqn:WPWS:PhiPrime}
  \ph'(x)=\alpha+\sum_{\xi-\eul\eta\in\Lambda}\frac{\eta}{(x-\xi)^2+\eta^2}.
 \end{equation}

Let $\Nev^+_h$ be the class of functions of bounded type on $\cpl^+$ whose mean type does not exceed $h$
 (see \cite{deBranges}, p.19 and p.26 for the definitions).
From an arbitrary Hermite--Biehler function $E$, one may build a prototypical example of de Branges space, namely \cite[prob.50]{deBranges}
 \begin{equation*}
  \Hil(E)=\Ens{f \mbox{ entire}}{\nm{\kfrac{f}{E}}_2<\infty\mbox{ and } f/E,\,f^*/E\in\Nev^+_0}.
 \end{equation*}
It is equipped with the norm $\nm{f}_{\Hil(E)}=\nm{\kfrac{f}{E}}_2$.
In fact, a theorem of de Branges \cite[th.23]{deBranges} shows that \emph{every de Branges space is isometrically equal to a space of the form
 $\Hil(E)$}, where $E$ is not unique in general.%
\footnote{%
The above definition of $\Hil(E)$ is the original one. 
Denoting by $H^2(\cpl^+)$ the Hardy space on the upper half-plane, one may prefer the well-known characterization
 $$\Hil(E)=\ens{f \mbox{ entire}}{\nm{\kfrac{f}{E}}_2<\infty\mbox{ and } f/E,\,f^*/E\in H^2(\cpl^+)}.$$}
 
The reproducing kernel in $\Hil(E)$ is given by \cite[th.19]{deBranges}
 \begin{equation}\label{Eqn:WPWS:NoyauReproduisant}
  k_\zeta(z)=\frac{E^*(z)\conj{E^*(\zeta)}-E(z)\conj{E(\zeta)}}{2\pi\eul(z-\conj{\zeta})}.
 \end{equation}
Therefore, if $E$ does not have real zeroes, then for $x,\xi\in\reel$
 $$k_\xi(x)=\frac{|E(x)|\,|E(\xi)|}{\pi}\,\frac{\sin(\ph(x)-\ph(\xi))}{x-\xi}.$$
In particular,
 \begin{equation}\label{Eqn:WPWS:Majorant}
  M(x)=\sqrt{k_x(x)}=\Inv{\sqrt{\pi}}\sqrt{\ph'(x)}|E(x)|
 \end{equation}
for all $x\in\reel$.\\

As mentioned above, every de Branges space admits a representation $\Hil(E)$ for an Hermite--Biehler function $E$.
Wlog $\alpha=0$ in the factorization of $E$.
Indeed, if $\alpha>0$, $E(z)$ may be replaced with $E_0(z)=E(z)\exi{\alpha z}\sin(\alpha(z+\eul))$, which is also an Hermite--Biehler function
 (due to its factorization).
Since $|E_0(z)|\simeq |E(z)|$ for $\Im z\geq 0$, $\Hil(E_0)$ and $\Hil(E)$ are equal with equivalent norms.
If in addition $\Hil(E)$ is a weighted Paley--Wiener space, then $E$ does not have real zeroes.
Otherwise \cite[prob.44]{deBranges} the relation $E(x)=0$ for an $x\in\reel$ would imply $f(x)=0$ for all $f\in\Hil(E)$, that is, $M(x)=0$, a contradiction.
For the above reasons, in this text the \emph{considered Hermite-Biehler functions} are those of the form~(\ref{Eqn:WPWS:FactoHermiteBiehler}) with $\alpha=0$, $m=0$, and
 $\Lambda\subset\cpl^-$.
In particular, \emph{every weighted Paley-Wiener space admits a representation $\Hil(E)$ for a considered HB-function $E$}.

\begin{exe}
 Given a considered HB-function of phase $\ph$, if $\ph'(x)\simeq 1$, then $\Hil(E)$ is obviously a weighted PW-space.
 The converse statement however does not hold in general, as we shall see later.
\end{exe}


\section{A prototypical example}\label{Sec:APE}

Lyubarskii and Seip's work \cite{LyuSeip} makes a bridge between Hermite--Biehler functions and a certain kind of potentials,%
\footnote{%
See \cite{Ransford} for potential theory of compactly supported measures on $\cpl$.}
 namely, potentials of measures
 of the form $m(x)\dif{x}$ for $m(x)$ measurable, positive, and $\simeq 1$.
Such a potential cannot be defined as $\int_{-\infty}^\infty\log|1-z/t|m(t)\dif{t}$, since this last integral does not exist.
This is due to the dominating term in the expansion
 $$\log|1-z/t|=-x/t-\sum_{n=2}^\infty(1/n)\Re(z^n)/t^n$$
for $|t|$ large.
It suggests to define
 \begin{equation*}
  \omega_m(z)=\int_{-\infty}^\infty\log^*|1-z/t|m(t)\dif{t}
 \end{equation*}
with $\log^*|1-z/t|=\log|1-z/t|+\cutoff(t)x/t$, where $\cutoff(t)=1-\chi_{[-1,1]}(t)$.

We first show that $\omega_m(z)$ is well-defined, indeed, that the above integral is absolutely convergent.
The previous expansion gives, for $|t|$ large,
 \begin{eqnarray}
  |\log^*|1-z/t|\,| &\leq& \sum_{n=2}^\infty (1/n)|z/t|^n\nonumber\\
                      &\leq& |z/t|^2\left(1/2+\sum_{n=1}^\infty (1/n)|z/t|^n\right)\nonumber\\
                      &  = & |z/t|^2(1/2-\log(1-|z/t|)\,).\label{Ineq:APE:LogStar}
 \end{eqnarray}
Since $m(x)\simeq 1$, it suffices to show that $\int_R^\infty -(1/t^2)\log(1-|z|/t)\dif{t}<\infty$ for $R$ large.
This last relation follows from the substitution $u=1-|z|/t$.
Therefore, $\int_{-\infty}^\infty|\log^*|1-z/t|\,|\dif{t}<\infty$.

The inequality (\ref{Ineq:APE:LogStar}) and the dominated convergence theorem also yield that $\omega_m$ is continuous.%
\footnote{%
Explicitly, let $C$ be the maximum of $|z|$ during the limiting process.
Since $m\simeq 1$, for $|t|>C$ a dominator may be obtained from
 $$|\log^*|1-z/t|\,|\leq (C^2/t^2)(1/2-\log(1-C/|t|)\,).$$
Regarding $\int_{-C}^C\log^*|1-z/t|m(t)\dif{t}$, it suffices to consider $\int_{-C}^C\log|t-z|m(t)\dif{t}$, that is,
 $\int_{-C-x}^{C-x}\log|t-\eul y|m(t+x)\dif{t}$.
The last integrand is dominated by $\sup m$ times the integrable function $\chi_{[-2C,2C]}(t)\,(\,\log^+(|t|+C)+\log^-|t|)$.}

Observe that, for $z\neq t$ and any appropriate branch of the logarithm,
 $$\partial_y\log^*|1-z/t|\,=\,\Re\partial_y\log(1-z/t)\,=\,\frac{y}{(x-t)^2+y^2}.$$
For $z\notin \reel$, the dominated convergence theorem then implies
 \begin{equation}\label{Eqn:APE:DeriveePotentiel}
  \partial_y\omega_m(z)=\pi P_m(z),
 \end{equation}
where $$P_m(z)=\KInv{\pi}\int_{-\infty}^\infty \frac{y}{(x-t)^2+y^2}m(t)\dif{t}$$ is the Poisson transform of $m(x)\dif{x}$.
In particular%
\footnote{%
In the above computation a dominator is easily obtained from the mean value theorem, for instance
 $MA/((x-t)^2+a^2)$, where $a=|y|/2$, $A=2|y|$, and $M=\sup m$.},
$|\partial_y\omega_m(z)|\simeq 1$.

Similarly, for $z\neq t$ and any appropriate branch of the logarithm,
 $$\partial_x\log^*|1-z/t|\,=\,\Re\partial_x\log(1-z/t)+\chi(t)/t\,=\,\frac{x-t}{(x-t)^2+y^2}+\frac{\chi(t)}{t}.$$
For $z\notin\reel$, the dominated convergence theorem implies%
\footnote{%
Again a dominator may easily be derived from the mean value theorem, for instance
 \begin{equation*}
  \left\{
   \begin{array}{lll}
    M\Kfrac{A|t|+A^2+y^2}{|t|\,((|t|-A)^2+y^2)} & \mbox{ if} & |t|\geq A\vspace{2mm}\\
    M+2M(A/y^2) & \mbox{ if} & |t|<A,
   \end{array}
  \right.
 \end{equation*}
where $A=|x|+1$ and $M=\sup m$.}
 \begin{equation*}
  \partial_x\omega_m(z)=\int_{-\infty}^\infty\left(\frac{x-t}{(x-t)^2+y^2}+\frac{\chi(t)}{t}\right)m(t)\dif{t}.
 \end{equation*}

Let us check%
\footnote{%
 The proof is mimicked on \cite[th.3.7.4]{Ransford}}
that $\Delta\omega_m=2\pi m(x)\dif{x}\dif{\delta_0(y)}$ in the sense of distribution, where $\delta_0$ denotes the $1$-dimensional
 Dirac measure at $0$.
Let $C^\infty_c$ be the space of smooth, compactly supported functions on $\cpl$.
Then, $\omega_m$ is identified with $C^\infty_c\to\cpl,\ \omega_m(\psi)=\int_\cpl\omega_m(z)\psi(z)\dif{A}(z)$, where $\dif{A}$
 is the element of area.
Observe that $\omega_m(\psi)$ well-defined, since $\omega_m$ is continuous.
By the definition of the distributional Laplacian,
 $$\Delta\omega_m(\psi)=\omega_m(\Delta\psi)=\int_\cpl\omega_m(z)\Delta\psi(z)\dif{A}(z),$$
and hence the harmonicity of $\omega_m$ in $\cpl\setminus\reel$ implies
 $$\Delta\omega_m(\psi)=\lim_{\eps\downarrow 0}\int_{|\Im z|>\eps}\omega_m\Delta\psi-\psi\Delta\omega_m\dif{A}.$$
This last integral may be taken on $\ens{z\in\cpl}{|z|<R\mbox{ and }|\Im z|>\eps}$ for $R$ large enough, since $\psi$ is compactly supported.
Consequently, Green's identity and the equation~(\ref{Eqn:APE:DeriveePotentiel}) imply that $\Delta\omega_m(\psi)$ is equal to
 $$\lim_{\eps\downarrow 0}\left(-\int_{-\infty}^\infty\omega_m(x+\eul\eps)\dpar{y}\psi(x+\eul\eps)\dif{x}+\int_{-\infty}^\infty \psi(x+\eul\eps)\pi P_m(x+\eul\eps)\dif{x}\ +\right.$$
 $$\left.+\int_{-\infty}^\infty\omega_m(x-\eul\eps)\dpar{y}\psi(x-\eul\eps)\dif{x}-\int_{-\infty}^\infty\psi(x-\eul\eps)\pi P_m(x-\eul\eps)\dif{x}\right).$$
Observe the first and third integrands in the last expression are continuous and compactly supported.
The dominated convergence theorem and the relation $P_m(x+\eul\eps)=-P_m(x-\eul\eps)$ then imply
 $$\Delta\omega_m(\psi)=\lim_{\eps\downarrow 0}\int_{-\infty}^\infty\pi(\,\psi(x+\eul\eps)+\psi(x-\eul\eps)\,)P_m(x+\eul\eps)\dif{x}.$$
Since the above integrand is bounded (independent of $\eps$) and compactly supported, the dominated convergence theorem yields
 \begin{equation}\label{Eqn:APE:DeltaPotentiel}
  \Delta\omega_m(\psi)=\int_{-\infty}^\infty2\pi\psi(x)m(x)\dif{x}=\int_\cpl2\pi\psi(z)m(x)\dif{x}\dif{\delta_0(y)},
 \end{equation}
as claimed.%
\footnote{%
Let us mention that our previous considerations may easily be generalized to any measurable function $m$ satisfying $|m|\lesssim 1$.
\comm{verifier soigneusement}
In fact, $\omega_m$ is then absolutely convergent, since its integrand is $\lesssim |\log^*|1-z/t|\,|$.}

\begin{exe}
 For $m(x)=1$ and $z\notin\reel$, $\partial_y\omega_1(z)=\pi\,\sgn y$, while $\partial_x\omega_1(z)=0$.
 Hence, $\omega_1(z)=\pi|y|+C$.
 By continuity, this last relation applies for all $z\in\cpl$.
 Since $\omega_1(0)=0$, we deduce
  \begin{equation*}
   \omega_1(z)=\pi|y|.
  \end{equation*}
 It is easily seen that $\Delta\omega_1=2\pi\dif{\delta_0}(y)\dif{x}$ in the sense of distribution.
\end{exe}

\kTroisEtoiles
Lyubarskii and Seip \cite{LyuSeip} made a bridge between $\omega_m(z)$ and an Hermite--Biehler function by counterbalancing $m(t)\dif{t}$ with a sum of masses.

Let  $\cdots<x_{-1}<x_0<x_1<\cdots$ be the partition of $\reel$ defined by
 $x_0=0$ and $\int_{x_k}^{x_{k+1}}m(t)\dif{t}=1$.
Observe that $x_{k+1}-x_k\simeq 1$, since $m(t)\simeq 1$.
Denoting by $\dif{\delta_\xi}$ the Dirac measure at $\xi\in\reel$, consider the auxiliary measure
 $$\dif{\mu}(t)=m(t)\dif{t}-\sum_{k=-\infty}^{\infty}\dif{\delta_{\xi_k}}(t),$$
where $\xi_k\in(x_k,x_{k+1})$ will be chosen later.

Let $f_\mu(x)=\int_{0^+}^x\dif{\mu}(t)$ (including the endpoint $x$) and let $F_\mu(x)=\int_0^x f_\mu(t)\dif{t}$.
For $a<b$ and $z\notin\reel$, two integrations by parts give%
\footnote{%
The formula below may also be derived by applying twice Fubini's theorem to
 $$\int_a^b\int_0^t\int_{0^+}^s\partial_t^2\log|1-z/t|\dif{\mu}(r)\dif{s}\dif{t},$$
and hence is justified by the finiteness of $\int_a^b\int_0^t\int_{0^+}^s|\:\partial_t^2\log|1-z/t|\:|\dif{|\mu|}(r)\dif{s}\dif{t}$.}
 \begin{equation*}
  \int_{a^+}^b\log|1-z/t|\dif{\mu}(t) = \int_a^b F_\mu(t)\partial_t^2\log|1-z/t|\dif{t} + R(a,b),
 \end{equation*}
where $R(a,b)=(f_\mu(t)\log|1-z/t|-F_\mu(t)\partial_t\log|1-z/t|)\ |_a^b$.

Since $\xi_k\in(x_{k},x_{k+1})$, clearly $f_\mu(x_k)=0$ for all $k$.
It follows that $f_\mu$ is a bounded function on $\reel$.
Moreover, we may choose $\xi_k$ so that $F_\mu(x_k)=0$ for all $k\in\rel$.
This is achieved when $\int_{x_k}^{\xi_k}\int_{x_k}^t m(s)\dif{s}\dif{t}+\int_{\xi_k}^{x_{k+1}}(\int_{x_k}^t m(s)\dif{s}-1)\dif{t}=0$,
 that is, when
  $$\xi_k=\int_{x_k}^{x_{k+1}}tm(t)\dif{t}.$$
Doing so, $\xi_k$ belongs to $(x_k,x_{k+1})$, is bounded away from $\sgl{x_k,x_{k+1}}$, and hence $\xi_{k+1}-\xi_k\simeq 1$.
In addition, $F_\mu$ is bounded on $\reel$.

On the one hand, we deduce that $R(a,b)\to 0$ when $a\to-\infty$ and $b\to\infty$.
On the other hand, using any appropriate branch of the logarithm,
 \begin{eqnarray*} \int_a^b F_\mu(t)\partial_t^2\log|1-z/t|\dif{t}
  &=& \int_a^b F_\mu(t)\,\Re\partial_t^2\log(1-z/t)\dif{t}\\
  &=& -\int_a^b F_\mu(t)\Re \Inv{(t-z)^2}\dif{t}+\int_a^b\frac{F_\mu(t)}{t^2}\dif{t}.
 \end{eqnarray*}
Since $F_\mu$ is bounded and since $F_\mu(t)\simeq t^2$ in a neighborhood of $x_0=0$, we conclude that
 $\int_{-\infty}^\infty\log|1-z/t|\dif{\mu}(t)$
is well-defined (as an improper integral) and satisfies
 \begin{equation*}
  \left|\int_{-\infty}^\infty\log|1-z/t|\dif{\mu}(t)\right|\lesssim 1
 \end{equation*}
when $\Im z\gg 0$ (i.e., when $\Im z>\eps>0$ for a certain $\eps>0$).

Integration by parts%
\footnote{%
It is justified by the finiteness of $\int_a^b\sgn{t}\int_{0^+}^t\chi(t)/t^2\dif{|\mu|(s)}\dif{t}$ for any $a<b$.}
 also gives
 $$\int_{-\infty}^\infty\frac{\chi(t)}{t}\dif{\mu(t)}=-f_\mu(-1)-f_\mu(1)+\int_{|t|>1}\frac{f_\mu(t)}{t^2}\dif{t}.$$

Since $f_\mu$ is bounded, the above integral is just a real constant $C$.
In total,
 \begin{equation*}
  \left|\int_{-\infty}^\infty\log^*|1-z/t|\dif{\mu}(t)-Cx\right|\lesssim 1
 \end{equation*}
when $\Im z\gg 0$.
Letting $\alpha=C-\sum_{|\xi_k|\leq 1}1/\xi_k$, it follows that
 $$|\,\omega_m(z)-\alpha x-\textstyle\sum_k\,(\log|1-z/\xi_k|+x/\xi_k)\,|\lesssim 1.$$
In other words, for $F_m(z)=\ex{\alpha z}\prod_k(1-z/\xi_k)\ex{z/\xi_k}$,
 $$|F_m(z)|\simeq\ex{\omega_m(z)}$$
when $\Im z\gg 0$.

Incidentally observe that $\alpha=\Klim_{R\to\infty}\int_{-R}^R\chi(t)m(t)/t\,\dif{t}-\sum_{|\xi_k|<R}1/\xi_k$.
This yields the factorization
 $$F_m(z)=\lim_{R\to\infty}\ex{z\int_{-R}^R\chi(t)m(t)/t\,\dif{t}}\prod_{|\xi_k|<R}(1-z/\xi_k).$$

Obviously, $|F_m(\bar{z})|=|F_m(z)|$, preventing $F_m$ to be in the Hermite--Biehler class.
Since in addition we need a relation valid on the whole upper half-plane, we consider
 instead $E_m(z)=F_m(z+\eul)$ and observe that for $z\in\cpl^+$
 $$|E_m(\bar{z})|=|F_m(z-\eul)|<|F_m(z+\eul)|=|E_m(z)|.$$
For all $z\in\cpl$, the mean value theorem (applied twice if $-1<\Im z<0$) implies $|\omega_m(z+\eul)-\omega_m(z)|\lesssim 1$,
 in other words $\ex{\omega_m(z+\eul)}\simeq\ex{\omega_m(z)}$.
Therefore,
 $$|E_m(z)|\simeq \ex{\omega_m(z)}\ \mbox{ when }\Im z\gg-1.$$

In total, we have proven the following \emph{multiplier lemma} \cite[lem.10]{LyuSeip}:

\begin{prop}\label{Prop:APE:HBCanonique}
 Let $m(x)\simeq 1$ be a measurable function.
 There exists an Hermite--Biehler function $E_m(z)$ satisfying
  $$|E_m(z)|\simeq\ex{\omega_m(z)}\ \mbox{ when }\Im z\gg-1.$$
 The zeroes of $E_m$ are simple and of the form $\xi_k-\eul$, where $\xi_{k+1}-\xi_k\simeq 1$.
 Here,
  $$E_m(z)=E_m(0)\lim_{R\to\infty}\ex{z\int_{-R}^R\chi(t)m(t)/t\,\dif{t}}\prod_{|\xi_k|<R}(1-z/(\xi_k-\eul)).$$
 In addition, $E_m(\bar{z})=E_m(z-2\eul)$ for all $z\in\cpl$.
 Finally, $$\ex{-\omega_m(z)-\eul\cjh{\omega}_m(z)}E_m(z)$$ is bounded and analytic on both $\cpl^+$ and $\cpl^-$.
\end{prop}

\TroisEtoiles
In the above context let $\ph$ be the phase of $E_m$.
By the equation (\ref{Eqn:WPWS:PhiPrime}),
 $$\ph'(x)=\sum_k\Inv{(x-\xi_k)^2+1}.$$
The condition $\xi_{k+1}-\xi_k\simeq 1$ then implies $\ph'(x)\simeq 1$.
Therefore $\Hil(E_m)$ is a weighted PW-space.
By the multiplier lemma, $E_m$ may be replaced with $\ex{\omega_m+\eul\cjh{\omega}_m}$ in the definition of $\Hil(E_m)$.
It follows that $\Hil(E_m)$ is equal with equivalent norms to the following space,
 $$PW(m)=\ens{f\mbox{ entire}}{\nm{f\ex{-\omega_m}}_2<\infty,\,f\ex{-\omega_m-\eul\tilde{\omega}_m},\,f^*\ex{-\omega_m-\eul\tilde{\omega}_m}\in\Nev^+_0},$$
equipped with the norm $\nm{f}_{PW(m)}=\nm{f\ex{-\omega_m}}_2$.
Since $\omega_m(\conj{z})=\omega_m(z)$, the above may be abbreviated as follows:
 denoting by $\Nev^\pm_h$ be the class of functions of bounded type on $\cpl\setminus\reel$ whose mean types on $\cpl^+$ and $\cpl^-$ do not exceed~$h$,
 $$PW(m)=\ens{f\mbox{ entire}}{\nm{f\ex{-\omega_m}}_2<\infty\mbox{ and }f\ex{-\omega_m-\eul\tilde{\omega}_m}\in\Nev^\pm_0}.$$

Let $L^2_h$ be the classical Paley-Wiener space, which consists of all square summable entire functions of exponential type at most $h$:
 $$L^2_h=\ens{f\mbox{ entire}}{\nm{f}_2<\infty\mbox{ and }|f(z)|\leq C_\eps\ex{(h+\eps)|z|}}.$$
By a theorem of Krein \cite[prob.37]{deBranges}, 
 $$L^2_h=\ens{f\mbox{ entire}}{\nm{f}_2<\infty\mbox{ and } f\in\Nev^\pm_h}.$$
Let $L^2_h[\Sigma]$ be the annihilator in $L^2_h$ of a given $\Sigma\subseteq\cpl$.
Clearly, it is a linear subspace of $L^2_h$, but not necessarily a de Branges space.

\begin{prop}\label{Prop:APE:Lifting}
 Given a measurable $m\simeq1$, assume $\tau>\sup m$.
 Let $E_{\tau-m}$ be given by Proposition~\ref{Prop:APE:HBCanonique} and $\Sigma$ be its zero set.
 Then, $f\mapsto fE_{\tau-m}$ is a bijection from $PW(m)$ to $L^2_{\pi\tau}[\Sigma]$ with equivalence of norms.
\end{prop}

\begin{pr}
 Proposition~\ref{Prop:APE:HBCanonique} implies that $|E_{\tau-m}(x)|\simeq\ex{\omega_{\tau-m}(x)}=\ex{-\omega_m(x)}$
  for $x\in\reel$.
 Thus $\nm{f\ex{-\omega_m}}_2\simeq\nm{fE_{\tau-m}}_2$ for any $f\in PW(m)$.
 Moreover, Proposition~\ref{Prop:APE:HBCanonique} and the Hermite--Biehler property imply that
  $|E_{\tau-m}|\ex{-\omega_{\tau-m}}$ is the modulus of an analytic, bounded function in $\cpl\setminus\reel$.
 In addition $\omega_\tau(z)=\pi\tau|\Im z|$.
 Thus, 
  \begin{eqnarray*}
   f\ex{-\omega_m-\eul\cjh{\omega}_m}\in\Nev_0^\pm 
    &\Leftrightarrow& f E_{\tau-m}\ex{\sgn{\Im z}\eul\pi\tau z}\in\Nev_0^\pm\\
    &\Leftrightarrow& f E_{\tau-m}\in\Nev_{\pi\tau}^\pm.
  \end{eqnarray*}
 In total, $f\in PW(m)$ if and only if $f E_{\tau-m}\in L^2_{\pi\tau}[\Sigma]$.
 Since the zeroes of $E_{\tau-m}$ are simple, it follows that $f\mapsto E_{\tau-m}f$ satisfies the desired properties.
\end{pr}

The above definition of $PW(m)$ may be simplified thanks to the previous proposition.

\begin{kcor}
  \rm $$PW(m)=\ens{f\mbox{ entire}}{\nm{f\ex{-\omega_m}}_2<\infty,\, |f(z)|\ex{-\omega_m(z)}\leq C_\eps\ex{\eps|z|}}.$$
\end{kcor}

\begin{exe}
 We have seen that $\omega_1(z)=\pi|y|$.
 Consequently,
  $$PW(1)=\ens{f\mbox{ entire}}{\nm{f}_2<\infty\mbox{ and }|f(z)|\leq C_\eps\ex{\pi|y|}\ex{\eps|z|}}.$$
 By a Phragm\'en-Lindel\"of argument \cite[lec.6 th.3]{Levin}, the above is the classical Paley--Wiener space $L^2_\pi$.
 It is indeed straightforward that $E_1(z)=(\kfrac{4}{\pi^2})\cos(\pi(z+\eul))$, and hence $|E_1(z)|\simeq|\exni{\pi z}|$
  on the closed upper half-plane, yielding $PW(1)=\Hil(\exni{\pi z})$, as expected.
\end{exe}

Given any real-entire $g$, $e^gPW(m)$ is also a weighted PW-space, equal to $\Hil(e^gE_m)$.
Since $\ex{g}E_m$ is an HB-function sharing its phase with $E_m$, the majorant-weight of $\ex{g}PW(m)$ is comparable with $|\ex{g(x)}E_m(x)|\simeq \ex{g(x)}\ex{\omega_m(x)}$.
Spaces of this last form may be regarded as a prototypical examples of weighted PW-space, as we shall see later.\\

\begin{rem}
 What about the majorant $M(z)$ of $PW(m)$ for $z$ in the upper half-plane?
 The equation~(\ref{Eqn:WPWS:NoyauReproduisant}) and the relation $E_m(\bar{z})=E_m(z-2\eul)$ imply, for $\Im z\geq 2$,
  $$M(z)^2=\Inv{4\pi\Im z}(|E_m(z)|^2-|E_m(z-2\eul)|^2)\simeq \frac{\ex{2\omega_m(z)}}{4\pi\Im z}(1-\ex{2\omega_m(z-2\eul)-2\omega_m(z)}).$$
 Furthermore, the mean value theorem and the relation $m\simeq 1$ imply
  $$\omega_m(z)-\omega_m(z-2\eul)\geq 2\inf_{\Im z-2\leq y\leq\Im z}y\int_{-\infty}^\infty \frac{m(t)}{(\Re z-t)^2+y^2}\dif{t}\gtrsim 1.$$
 Therefore,
  $M(z)\simeq\kfrac{\ex{\omega_m(z)}}{\sqrt{\Im z}}$
 when $\Im z\geq 2$.
 This last limitation is due to our choice to plot the zeroes of $E_m$ on the axis $\Im z=-1$.
 Clearly, by plotting them on $\Im z=-\eps$ for an arbitrary $\eps>0$, the above relation extends to $\Im z\gg 0$.
 In conclusion, the majorant of $\ex{g}PW(m)$ simultaneously satisfies $M(x)\simeq\ex{g(x)}{\ex{\omega_m(x)}}$ for $x\in\reel$, and
  \begin{equation}\label{Eqn:APE:MajorantDemiPlan}
   M(z)\simeq\frac{\ex{\Re g(z)}\ex{\omega_m(z)}}{\sqrt{\Im z}}\mbox{ for $\Im z\gg 0$}.
  \end{equation}
\end{rem}

\section{The mountain chain axioms} \label{Sec:MCA}

Let $\Hil=\Hil(E)$ be a weighted Paley--Wiener space, where $E$ is a considered HB-function.
In order to understand the structure of $\Hil$, Lyubarskii and Seip \cite{LyuSeip} studied the distribution of the zeroes of $E$.
These lasts determine the phase via the equation~(\ref{Eqn:WPWS:PhiPrime}).
Their work is entirely achieved by means of the following relation applied to several test functions: by the definition of weighted PW-space and the
 equation~(\ref{Eqn:WPWS:Majorant}),
 \begin{equation*}
  \int_{-\infty}^\infty\frac{|f(t)|^2}{|E(t)|^2}\dif{t}\simeq\int_{-\infty}^\infty\frac{|f(t)|^2}{|E(t)|^2\ph'(t)}\dif{t}
 \end{equation*}
for all $f\in\Hil$.

As a first instance, for $f=k_x/E(x)$, where $x\in\reel$ is arbitrarily fixed, the above gives
 \begin{equation*}
  \ph'(x)\simeq\int_{-\infty}^\infty\frac{\sin^2(\ph(t)-\ph(x))}{(t-x)^2\ph'(t)}\dif{t}.
 \end{equation*}
In particular, $E$ has infinitely many zeroes: otherwise, the equation (\ref{Eqn:WPWS:PhiPrime}) would imply $\ph'(t)\simeq 1/t^2$ for $|t|$ large,
 contradicting (for almost all $x$) the finiteness of previous integral.
Therefore, test functions of the form $f(z)=E(z)/(z-\lambda)$, where $\lambda$ is a zero of $E$, are available in $\Hil(E)$.
This yields
 \begin{equation*}
  -\Inv{\Im\lambda}\simeq\int_{-\infty}^\infty\Inv{|t-\lambda|^2\ph'(t)}\dif{t}.
 \end{equation*}
Tests functions of the form $f(z)=E(z)/(z-\lambda)(z-\lambda')$, where $\lambda$ and $\lambda'$ are zeroes of $E$, are also used in the original paper \cite[p.987]{LyuSeip}.
A succession of lemmas regarding the zeroes of $E$ are then deduced from the three resulting relations.
Their proofs are difficult, while involving only elementary tools.
They are not be presented here.
The purpose of the present section is to rephrase results to be used later and list them as axioms.\\

In the sequel $\Lambda$ denotes the zero sequence (sometimes, the zero set) of $E$, $\str_\delta=\sgl{-\delta<\Im z<0}$, and
 $\lambda(x)=\xi(x)-\eul\eta(x)$ denotes the element of $\Lambda$ closest to a given $x\in\reel$ (if there are several candidates
 for $\lambda(x)$, we select the one with the smallest $x$-coordinate.)
The lemma~1 in the original paper \cite{LyuSeip} (presented as the \emph{chief result}) states that, for $\delta>0$ sufficiently small,
\begin{enumerate}
 \item \label{Enu:MCA:ChiefSeparation}
  If $\lambda\in\Lambda\cap\str_\delta$, then $\dist{\lambda}{\Lambda\setminus\sgl{\lambda}}\simeq 1$;
 \item \label{Enu:MCA:ChiefPoisson}
  $\displaystyle\ph'(x)\simeq\frac{\min(1,\eta(x))}{\min(1,|x-\lambda(x)|^2)}.$
\end{enumerate}
Let us rephrase the latter conclusion.
For a given $x\in\reel$ and any $\delta\in(0,1]$, let us divide $\cpl^-$ into three regions: the half-disk $\I_\delta=\cpl^-\cap\sgl{|z-x|<\delta}$, the rest of the strip
 $\II_\delta=\str_\delta\setminus\I_\delta$, and the rest of the half-plane $\III_\delta=\cpl^-\setminus\str_\delta$.
The conclusion \ref{Enu:MCA:ChiefPoisson} is then equivalent to
 \begin{equation*}
  \ph'(x)\simeq
   \left\{
    \begin{array}{lll}
     \eta(x)/|x-\lambda(x)|^2 &\mbox{ if }& \lambda(x)\in\I_1,\\
     \eta(x) &\mbox{ if }& \lambda(x)\in\II_1,\\
     1 &\mbox{ if }& \lambda(x)\in\III_1.
    \end{array}
   \right.
 \end{equation*}
We may replace the index $1$ by $\delta$ in the above.
Because, $|x-\lambda(x)|\simeq 1$ if $\lambda(x)\in\II_\delta\cap\I_1$, and hence $\ph'(x)\simeq\eta(x)$ if $\lambda(x)\in\II_\delta$.
Similarly, both $\eta(x)$ and $|x-\lambda(x)|$ are $\simeq 1$ if $\lambda(x)\in\III_\delta\cap\I_1$, while $\eta(x)\simeq 1$ if $\lambda(x)\in\III_\delta\cap\II_1$,
 and hence $\ph'(x)\simeq 1$ if $\lambda(x)\in\III_\delta$.

By choosing $\delta$ sufficiently small, if $\lambda(x)\in\II_\delta$, then $\ph'(x)\simeq\eta(x)$ becomes small enough for the original lemma~2 to apply.
Then $|\lambda(x)-x|\lesssim 1$, and hence $|\lambda(x)-x|\simeq 1$ by the definition of $\II_\delta$.
Therefore, for $\delta$ small enough
 \begin{equation} \label{Eqn:MCA:Phi'}
  \ph'(x)\simeq
   \left\{
    \begin{array}{ll}
     \eta(x)/|x-\lambda(x)|^2 &\mbox{ if }\lambda(x)\in\str_\delta,\\
     1 & \mbox{ otherwise.}
    \end{array}
   \right.
 \end{equation}
\Rtr

For $\lambda\in\Lambda$, let $I_\lambda=\ens{x\in\reel}{\lambda(x)=\lambda}$.
Notice that $I_\lambda$ is a half-closed interval, bounded or unbounded, empty or nonempty, whose left endpoint is discarded.
We show that, if $\delta>0$ is sufficiently small, then $|I_\lambda|\simeq 1$ when $\lambda$ varies in $\Lambda\cap\str_\delta$ (a conclusion vacuously true
 if this last set is empty).

If $\delta>0$ is sufficiently small, then $I_\lambda\neq\emptyset$ for $\lambda\in\Lambda\cap\str_\delta$.
Otherwise, there would exist a sequence of $\lambda\in\Lambda$ going to the real axis and such that $I_\lambda=\emptyset$.
Along this sequence, $|\lambda(\Re\lambda)-\lambda|\to 0$, while $\lambda(\Re\lambda)\neq\lambda$.
This contradicts the chief result \ref{Enu:MCA:ChiefSeparation}.

Observe that, if $\delta>0$ is sufficiently small, then $\Re\lambda\in\I_\lambda$ for $\lambda\in\Lambda\cap\str_\delta$.
Otherwise, there would be a sequence of $\lambda\in\Lambda$ going to the real axis and satisfying $\Re\lambda\notin\I_\lambda$, $I_\lambda\neq\emptyset$.
Let $x_\lambda$ be the endpoint of $I_\lambda$ closest to $\lambda$.
Then, the circle centered at $x_\lambda$ of radius $|\lambda-x_\lambda|$ would contain an element of $\Lambda$ different from $\lambda$, lying
  on the small arc joining $\lambda$ with the real axis (by the definition of $\I_\lambda$).
This contradicts the chief result \ref{Enu:MCA:ChiefSeparation}.

Consequently, if $\delta>0$ is sufficiently small, then $|I_\lambda|\gtrsim 1$ when $\lambda\in\Lambda\cap\str_\delta$.
Otherwise, there would exist a sequence of $\lambda\in\Lambda$ going to the real axis such that $|I_\lambda|\to 0$ and $\Re\lambda\in I_\lambda$.
In particular, the distance between an endpoint of $I_\lambda$ and $\lambda$ would go to zero.
But this last is also the distance between the same endpoint and another element of $\Lambda$, which goes closer and closer to $\lambda$.
Again, this contradicts the chief result \ref{Enu:MCA:ChiefSeparation}.

Finally, $|I_\lambda|\lesssim 1$ when $\lambda\in\Lambda\cap\str_\delta$.
Otherwise $|I_\lambda|\to\infty$ along a sequence of $\lambda\in\Lambda\cap\str_\delta$.
Then, the distance between a portion of $I_\lambda$ and $\lambda$ is going to infinity, and hence the restriction of $\ph'(x)$ to this last portion is going to zero, by the equation~(\ref{Eqn:MCA:Phi'}).
These two facts contradict the original lemma~2.
In total, the proof of our claim is complete.

Let us add that, if $\delta>0$ is sufficiently small, then $\Re\lambda\in I_\lambda$ is bounded away from $\partial I_\lambda$ when $\lambda$ varies in $\Lambda\cap\str_\delta$.
Otherwise, there would be a sequence of $\lambda\in\Lambda$ going to the real axis and satisfying $|\Re\lambda-x_\lambda|\to 0$, where $x_\lambda$ is
 the endpoint of $I_\lambda$ closest to $\lambda$.
In total, $|\lambda-x_\lambda|$ would go to zero.
Since the circle of radius $|\lambda-x_\lambda|$ centered at $x_\lambda$ contains an element of $\Lambda$ different from $\lambda$, this
 contradicts the chief result \ref{Enu:MCA:ChiefSeparation}.\\

Let us summarize our considerations as follows.
In the sequel, a piecewise continuous function $\reel\to\reel$ is called a \emph{mountain chain} if its graph consists of a succession of
 continuous pieces, called \emph{mountains} and \emph{plateaux}, satisfying the following conditions:
\begin{itemize}
 \item
  Each mountain has a Poissonian shape with two sides and a summit;
 \item
  The bases of the montains are $\simeq 1$;
 \item
  The summits have level more than 1; horizontally, they are bounded away from the endpoints of the mountain bases;
 \item
  The plateaux consist of horizontal segments of level $1$, without restriction on their lengths (finite or infinite).
\end{itemize}

The previous considerations yield:
\emph{%
 Let $\gamma_\delta(x)$ be the right-hand side of the equation (\ref{Eqn:MCA:Phi'}).
 If $\delta$ is sufficiently small, then $\gamma_\delta$ is a mountain chain.
}
In fact, $x\in\reel$ is under a mountain of summit $(\xi(x),1/\eta(x))$ if $\lambda(x)=\xi(x)-\eul\eta(x)$ is in the strip $\str_\delta$,
 $x$ is under a plateau otherwise.\\

In the case of weighted PW-spaces, since $\ph'\simeq\gamma_\delta$, the continuity of $\ph'$ gives a weak limitation on the growth of
 the summits of $\gamma_\delta$.

\begin{prop}\label{Prop:MCA:CroissanceSommets}
 Let $(\xi_1,1/\eta_1)$, $(\xi_2,1/\eta_2)$ be two consecutive summits of a mountain chain.
 Assume that this last is comparable with a continuous function.
 Then, $\eta_2/\eta_1\simeq 1$ independent of the choice of the consecutive summits.
\end{prop}

\begin{pr}
 Let $\gamma(x)$ and $f(x)$ be the mountain chain and the continous function in question, so $a\gamma(x)\leq f(x)\leq A\gamma(x)$ for some positive
  constants $a$ and $A$.

 In the case where the consecutive mountains are not separated by a plateau, let $x_0\in\reel$ be at their junction.
 Then, $a\gamma(x_0^+)\leq f(x_0)\leq A\gamma(x_0^-)$, while $a\gamma(x_0^-)\leq f(x_0)\leq A\gamma(x_0^+)$.
 Therefore, $\gamma(x_0^+)/\gamma(x_0^-)\simeq 1$, in other words
  $$\frac{\eta_2}{\eta_1}\left(\frac{\eta_1^2+(x_0-\xi_1)^2}{\eta_2^2+(x_0-\xi_2)^2}\right)\simeq 1.$$
 Since $\eta_j\leq 1$ and $|x_0-\xi_j|\simeq 1$ for $j=1,2$, the result follows.

 Consider a mountain of summit $(\xi,1/\eta)$ adjacent to a plateau, and let $x_0\in\reel$ be at their junction.
 Again, $\gamma(x_0^+)/\gamma(x_0^-)\simeq 1$, that is, $\gamma(x_0^\pm)\simeq 1$.
 Therefore, $\eta\simeq 1$.
 In particular, in the case where the consecutive mountains are separated by a plateau, $\eta_1\simeq 1\simeq \eta_2$, and the result follows.
\end{pr}
\begin{cor}
 Let $\sgl{(\xi_k,1/\eta_k)}$ be an indexation of the summits in a given mountain chain in increasing order of $x$-coordinates,
  where $-\infty\leq M\leq k\leq N\leq\infty$.
 Assume that the mountain chain is comparable with a continuous function.
 Then, $|\log\eta_k-\log\eta_l|\lesssim |k-l|.$
\end{cor}

A fortiori, $|\log\eta_k-\log\eta_l|\lesssim |\xi_k-\xi_l|$.
Lyubarskii and Seip obtained a better limitation on the growth of the summits for weighted PW-spaces \cite{LyuSeip}.
In this last case, by the original equation~(16), alternatively, by the original lemma~3,
 $|\log\eta_k-\log\eta_l|\lesssim \log|\xi_k-\xi_l|+1$ when $|\xi_k-\xi_l|$ is large enough.
In search for axioms, we retain that
 \begin{equation} \label{Eqn:MCA:CroissanceSommets1}
  \log\eta_k-\log\eta_l=O(|\xi_k-\xi_l|^{1-\eps})
 \end{equation}
uniformly in $l$ when $|k-l|\to\infty$, where $\eps>0$ is arbitrarily small.\\

Let $E$ be a considered HB-function of phase $\ph$.
In the sequel we shall say that $E$ satisfies the \emph{mountain chain axioms} if
 \emph{\begin{enumerate}
  \item\label{Axi:MCA:ChaineMontagnes}
   $\ph'$ is comparable with a mountain chain whose summits correspond to the zeroes of $E$ in a critical band, in the sense of the
   equation (\ref{Eqn:MCA:Phi'});
  \item
   These last zeroes are all simple;
  \item
   The summits of the mountain chain have a limited growth, in the sense of the equation (\ref{Eqn:MCA:CroissanceSommets1}).
 \end{enumerate}}

The following observation, though naive, is useful.

\begin{prop}\label{Prop:MCA:SourceMontagnes}
 Let $\xi(x)-\eul\eta(x)$ be the zero of $E$ closest to $x$.
 If $E$ satisfies the mountain chain axioms, then the derivative of its phase satisfies
  $$\ph'(x)=\frac{\eta(x)}{(x-\xi(x))^2+\eta(x)^2}+O(1)$$
 when $x$ varies below the mountains.
\end{prop}

\begin{pr}
 Recall that the summits of the mountains are at a horizontal distance $\gtrsim 1$ from their respective endpoints.
 In addition, the mountains have height at least $1$ and have Poissonian shape.
 It follows that their endpoints have height $\lesssim 1$.
 In particular, denoting by $(a_x,b_x)$ the base of the mountain over $x$, $\ph'(a_x^+)\lesssim 1$ and $\ph'(b_x^-)\lesssim 1$.

 Let us decompose the right-hand side of~(\ref{Eqn:WPWS:PhiPrime}) in
  $$\ph'(t)=\frac{\eta(x)}{(t-\xi(x))^2+\eta(x)^2}+L(t)+R(t)+U(t),$$
 where $L(t)$ is the sum of the terms with $\xi<a_x$, $R(t)$ is the sum of the terms with $\xi>b_x$, and $U(t)$ is the sum of the remaining terms.
 Clearly, $L(t)$ is decreasing on $(a_x,b_x)$, and hence $L(x)\leq\ph'(a_x^+)\lesssim 1$.
 Similarly $R(t)$ is increasing on $(a_x,b_x)$, and hence $R(x)\leq\ph'(b_x^-)\lesssim 1$.
 It remains to consider $U(x)$.
 Notice that $\xi(x)-\eul\eta(x)$ belongs to the critical strip $\str_\delta$, since $x$ is under a mountain.
 Consequently, the zeroes of $E$ contributing to $U(x)$ lie under $\str_\delta$.
 If $\lambda=\xi-\eul\eta$ is such a zero, let us show that its contribution $p(t)=\eta/((t-\xi)^2+\eta^2)$
  satisfies $p(t)\simeq p(a_x)$ on $[a_x,b_x]$ independent of $x$.
 For any $s,t\in[a_x,b_x]$, since $b_x-a_x\leq\ell$ for a certain $\ell$ independent of the mountain, $p(s)\geq\eta/(\ell^2+\eta^2)$, while $p(t)\leq 1/\eta$,
  and hence $p(t)/p(s)\leq (\ell^2/\eta^2)+1$.
 This proves our claim, since $\eta\geq\delta$.
 Therefore, $U(x)\simeq U(a_x)\leq\ph'(a_x)\lesssim 1$, and the result follows.
\end{pr}

\begin{sch}
 The proof shows more: the difference $$\ph'(x)-\Kfrac{\eta(x)}{(\xi(x)-x)^2+\eta(x)^2}$$ under a mountain is arbitrarily small if the mountain is
  sufficiently high.
\end{sch}

\begin{krem}
 Consider an HB-function $E$ satisfying the mountain chain axioms with a critical strip $\str_\delta$.
 Let $\ph$ be its phase and $\sgl{\xi_k-\eul\eta_k}$ be its zeroes.
 Observe that $|\xi_k|\gtrsim k$, and hence $\sum_{\eta_k<\delta}(\xi_k^2+\eta_k^2)^{-1}<\infty$.
 Moreover, the equation~(\ref{Eqn:WPWS:PhiPrime}) implies $\sum_{\eta_k\geq\delta}\delta(\xi_k^2+\eta_k^2)^{-1}\leq\ph'(0)<\infty$.
 The canonical product \cite[p.302]{Rudin} of $E$ thus has genus~1, and hence the factorization (\ref{Eqn:WPWS:FactoHermiteBiehler}) reduces to
  $$E(z)=C z^m \ex{h(z)}\exni{\alpha z}\prod_{\lambda\in\Lambda}(1-z/\lambda)\ex{z\Re(1/\lambda)}.$$
\end{krem}
\kTroisEtoiles
In the sequel we shall say that a de Branges space is an \emph{MC-space} if it has a representation $\Hil(E)$ for a considered HB-function satisfying the mountain chain axioms, in
 other words, if it is isometrically equal to $\Hil(E)$ for such an $E$.
In the first part of their paper \cite{LyuSeip} Lyubarskii and Seip proved that \emph{any weighted PW-space is an MC-space}, as explained above.
Then, they deduced results about MC-spaces, which they formulated in terms weighted PW-spaces.
However, it appears in Section~\ref{Sec:NSC} (corresponding to the original theorem~4) that there exists MC-spaces which are not weighted PW-spaces.
For this reason, in the following sections we shall reformulate Lyubarskii and Seip's theorems in terms of MC-spaces.
We shall elaborate their proofs in detail.

\section{Generic form of the majorants-weights}\label{Sec:GFMW}

Let $E$ be a considered HB-function of phase $\ph$ and suppose $\ph'(x)\simeq 1$ on the real line.
In this basic case, $\Hil(E)$ is a weighted PW-space by the equation~(\ref{Eqn:WPWS:Majorant}).
Let $H(z)=\log|E(z)|$.
Since the zeroes of $E$ lie on the lower half-plane, $H(z)$ is harmonic in an open neighborhood of $\adh{\cpl^+}$.
In particular, $\Delta H$ evaluated at $x+\eul|y|$ vanishes for all $x+\eul y\in\cpl$.

Now consider the composite function $H(x+\eul|y|)=\log|E(x+\eul|y|)|$.
Since $E$ is an Hermite--Biehler function, $H(x+\eul|y|)$ is subharmonic in $\cpl$.
Let us compute its distributional Laplacian.
Clearly, $\dpar{x}^2(H(x+\eul|y|)\,)=\dpar{x}^2 H(x+\eul|y|)$.
Moreover, $\dpar{y}(H(x+\eul|y|)\,)=\dpar{y}H(x+\eul|y|)\sgn y$, and hence
 $$\dpar{y}^2(H(x+\eul|y|)\,)=\dpar{y}^2H(x+\eul|y|)+2\dpar{y}H(x+\eul|y|)\dif{\delta_0(y)}\dif{x}.$$
In total,
 $$\Delta(H(x+\eul|y|))=\Delta H(x+\eul|y|)+2\dpar{y}H(x+\eul|y|)\dif{\delta_0(y)}\dif{x}=2\dpar{y}H(x)\dif{x}\dif{\delta_0(y)}.$$
Since $\log E(x)=H(x)-\eul\ph(x)$ has an analytic extension in a neighborhood of $\reel$, the Cauchy--Riemann equations yield $\dpar{y}H(x)=\ph'(x)$, and hence
 $$\Delta(H(x+\eul|y|)\,)=2\ph'(x)\dif{x}\dif{\delta_0(y)}.$$

Since $\ph'(x)\simeq 1$, $\omega_{\ph'/\pi}$ is well-defined.
The equation~(\ref{Eqn:APE:DeltaPotentiel}) and the above then provide a \emph{Riesz decomposition} \cite[p.991]{LyuSeip}
 $$H(x+\eul|y|)=\log|E(x+\eul|y|)|=h(z)+\omega_{\ph'/\pi}(z),$$
where $h$ is harmonic and $z=x+\eul y$ varies in $\cpl$.
In addition, $h(z)=h(\bar{z})$, that is, $h(z)=\Re g(z)$ for a real-entire function $g$.
We have proven:

\begin{prop}\label{Prop:GFMW:CasFacile}
 Let $E=|E|\ex{-\eul\ph}$ be a considered Hermite--Biehler function.
 Suppose $\ph'(x)\simeq 1$.
 Then, there exists a real-entire function $g(z)$ such that
  $$|E(x+\eul|y|)|=\ex{\Re g(z)}\ex{\omega_{\ph'/\pi}(z)}$$
 for all $z=x+\eul y$ in $\cpl$.
\end{prop}

Invoking the equation~(\ref{Eqn:WPWS:Majorant}), we deduce

\begin{cor}
 Let $M$ be the majorant of $\Hil(E)$.
 In the above circumstances, the restriction of $M$ to the real axis admits a representation
  $$M(x)\simeq\ex{g(x)}\ex{\omega_m(x)}$$
 for a certain $m(x)\simeq 1$.
\end{cor}

In fact, this last conclusion may be generalized to every MC-space, as we now show (with tears!) until the end of the section.

Let $E$ be an HB-function satisfying the mountain chain axioms.
We denote by $\ph$ its phase and by $\Lambda$ its family of zeroes (with multiplicity).
Let $\str_\delta$ be the strip in the equation~(\ref{Eqn:MCA:Phi'}), so the elements of $\Lambda\cap\str_\delta$ are simple zeroes lying under
 the summits of the mountains.
We suppose that $\ph'$ is not $\simeq 1$, so the zeroes of $E$ accumulate to the real axis.
In particular $\Lambda\cap\str_\delta$ contains infinitely many elements.
For $z\in\cpl$, let $\lambda_0(z)$ be the zero of $\Lambda\cap\str_\delta$ closest to $\Re z$ 
 (in case of equality, we may choose the one with the smallest $x$-coordinate).
Let $\sgl{\lambda_k(z)}$ be the resulting indexation of the zeroes of $\Lambda\cap\str_\delta$ in order of $x$-coordinates, where 
 $-\infty\leq M\leq k\leq N\leq\infty$. 
We write $\lambda_k(z)=\xi_k(z)-\eul\eta_k(z)$ and $z=x+\eul y$.
By the axiom~\ref{Axi:MCA:ChaineMontagnes} the zeroes in $\str_\delta$, which correspond to the summits, are at horizontal
 distance $\gtrsim 1$ one from another, and hence
 \begin{equation}\label{Eqn:MCA:DistEquidistribution}
  |x-\xi_k(z)|\gtrsim |k|.
 \end{equation}

Lyubarskii and Seip's strategy \cite{LyuSeip} is to shift down the zeroes in $\str_\delta$ for constructing an auxiliary HB-function whose phase has derivative $\simeq 1$.
We thus define $\tilde{\lambda}_k(z)=\xi_k(z)-\eul\delta$ and consider the auxiliary function
 $$F(z)=E(z)\left(\frac{1-z/\tilde{\lambda}_0(z)}{1-z/\lambda_0(z)}\right)\prod_{k\neq 0}\frac{1-z/\tilde{\lambda}_k(z)}{1-z/\lambda_k(z)}.$$
Let us show that the above infinite product is well-defined.
Without forgetting the dependence in $z$, its general factor is
 $1+\eul z\,\kfrac{(\delta-\eta_k)}{(\tilde{\lambda}_k(z-\lambda_k))}.$ 
The relation~(\ref{Eqn:MCA:DistEquidistribution}) implies that
 $\sum_k|\tilde{\lambda}_k(z-\lambda_k)|^{-1}$ is bounded when $z$ varies in a compact.
Consequently, the product is absolutely convergent on each compact.

Observe that for $z=x+\eul y$ in the closed upper half-plane,
 $$\left|\frac{\lambda_k(x)}{\tilde{\lambda}_k(x)}\right|
    = \left|\frac{\lambda_k(z)}{\tilde{\lambda}_k(z)}\right|
    \leq \left|\frac{1-z/\tilde{\lambda}_k(z)}{1-z/\lambda_k(z)}\right|
    \leq \left|\frac{\tilde{\lambda}_k(z)-z}{\lambda_k(z)-z}\right|.$$
On the one hand,
 $$\left|\frac{\lambda_k(x)}{\tilde{\lambda}_k(x)}\right|^2
    \geq 1-\frac{\delta^2}{|\tilde{\lambda}_k(x)|^2}
    \geq\left(1-\frac{\delta^2}{|\tilde{\lambda}_k(x)|^2}\right)^2.$$
The separation of the $\xi_k(x)$'s implies that $\prod_{k\neq 0}|\lambda_k(x)/\tilde{\lambda}_k(x)|$ is well-defined and $\gtrsim 1$.
On the other hand,
 $$\left|\frac{z-\tilde{\lambda}_k(z)}{z-\lambda_k(z)}\right|^2
    \leq 1+\frac{2y+1}{(x-\xi_k(z))^2+y^2}
    \leq \left(1+\frac{2y+1}{(x-\xi_k(z))^2+y^2}\right)^2.$$
The relation (\ref{Eqn:MCA:DistEquidistribution}) then implies the existence of a constant $C$ such that
 $$\left|\frac{z-\tilde{\lambda}_k(z)}{z-\lambda_k(z)}\right| \leq 1+\frac{2y}{C k^2+y^2}+\Inv{C k^2}.$$
Since $\Ksum_{k\neq 0}\Inv{k^2}$ is finite and $\Ksum_{k\neq 0}\Kfrac{y}{C k^2+y^2}$ is bounded when $y$ varies in $[0,\infty)$, we deduce that $\Kprod_{k\neq 0}|z-\tilde{\lambda}_k(z)|/|z-\lambda_k(z)|$ is well-defined and $\lesssim 1$.

In total, $|\kfrac{F(z)}{E(z)}|\simeq \kfrac{|z-\tilde{\lambda}_0(z)|}{|z-\lambda_0(z)|}$.
Assume $|z-\lambda_0(z)|\geq 1$. 
Then, 
 $$1\leq\frac{|z-\tilde{\lambda}_0(z)|^2}{|z-\lambda_0(z)|^2}
    \leq 2+\frac{2y}{|z-\lambda_0(z)|^2} \leq 2+2\min(y,y^{-1})\leq 4.$$
In general, we thus obtain 
 \begin{equation}\label{Eqn:GFMW:RatioHBFunctions}
  \frac{|F(z)|}{|E(z)|}\simeq\Inv{\min(1,|z-\lambda_0(z)|)}, 
 \end{equation}
in other words $|F(z)|/|E(z)|\simeq\kInv{\min(1,\dist{z}{\Lambda\cap\str_\delta})}$ for any $z$ in the closed upper-half plane.

The right-hand side in the last relation may be rewritten $\ex{\log^-|z-\lambda_0(z)|}$.
Let us show that it is $\lesssim \ex{-P_{\log\eta_0}(z)}$, where $P_{\log \eta_0}(z)$ is the Poisson transform of $\log\eta_0(t)$ at $z$.
This last is well-defined, since the assumption (\ref{Eqn:MCA:CroissanceSommets1}) implies $|\log\eta_0(t)|\lesssim|\xi_0(t)-\xi_0(0)|^{1-\eps}+1$, 
 where $|\xi_0(t)-\xi_0(0)|\leq 2|t-\xi_0(0)|$.
We want to prove that $\log^-|z-\lambda_0(z)|+P_{\log\eta_0}(z)$ is bounded above.
It suffices to show that $P_{\log\eta_0}(z)-\log\eta_0(z)$ is bounded for $|x-\xi_0(x)|\leq 1$ and $y\leq 1$, in other words that in these
 circumstances
 $$\frac{y}{\pi}\int_{-\infty}^\infty\frac{\log\eta_0(x+t)-\log\eta_0(x)}{t^2+y^2}\dif{t}\leq C<\infty.$$
The assumption (\ref{Eqn:MCA:CroissanceSommets1}) implies $|\log\eta_0(x+t)-\log\eta_0(x)|\lesssim|\xi_0(x+t)-\xi_0(x)|^{1-\eps}$, while the condition
 $|x-\xi_0(x)|\leq 1$ implies $|\xi_0(x+t)-\xi_0(x)|\leq 2|t|+2$.
The result then follows from the continuity of 
 $$(\kfrac{y}{\pi})\int_{-\infty}^\infty\kfrac{(2|t|+2)^{1-\eps}}{(t^2+y^2)}\dif{t}$$
on the compact set $[0,1]$.

We have proven that $1\lesssim |F(z)/E(z)|\lesssim \ex{-P_{\log \eta_0}(z)}$. 
Since the right-hand side of this last relation is the modulus of a function in $\Nev_0^+$, it follows that $F/E\in\Nev_0^+$, while $E/F\in\Nev_0^+$
 (see \cite{deBranges}, th.\@ 9 and 11).
In particular,

\begin{prop}\label{Prop:GFMW:InterchangementEF}
 For all functions $f$, $f/E\in\Nev_0^+$ if and only if $f/F\in\Nev_0^+$.
\end{prop}

It remains to analyse the relation (\ref{Eqn:GFMW:RatioHBFunctions}) on the real axis.
If $x$ is under a mountain, $\lambda_0(x)=\lambda(x)$, and hence the equation~(\ref{Eqn:MCA:Phi'}) implies
 $$|F(x)/E(x)|^2\simeq\ph'(x)/\eta(x).$$
If $x$ is under a plateau, $|x-\lambda_0(x)|\geq\delta$, and hence
 $|F(x)/E(x)|^2\simeq 1$, while $\ph'(x)\simeq 1$.
In both cases
 \begin{equation}\label{Eqn:GFMW:RatioFEAxeReel}
  |F(x)/E(x)|^2\simeq \ph'(x)/\sigma(x),
 \end{equation}
where $\sigma(x)=\min(1,\eta(x))$.

Due to its factorization $F$ is a considered HB-function.
Let $\psi$ be its phase.
Consider the development of $\psi'$ given by the equation~(\ref{Eqn:WPWS:PhiPrime}).
The proposition~\ref{Prop:MCA:SourceMontagnes} implies that the contribution to $\psi'(x)$ of the zeroes outside the critical strip is bounded
 when $x$ varies in $\reel$.
The contribution of the shifted zeroes is also bounded, due to the structure of the original mountain chain.
In addition, if $x$ is far from the mountains, $\psi'(x)\geq\ph'(x)\simeq 1$; otherwise, the proximity of a mountain suffices for $\psi'(x)\gtrsim 1$. 
Therefore, $\psi'(x)\simeq 1$ for all $x\in\reel$. 
By the previous considerations, $\Hil(F)$ is a weighted PW-space whose majorant-weight is comparable with $|F(x)|\simeq\ex{h(x)}\ex{\omega_{\psi'/\pi}(x)}$
 for a real-entire $h$.
Therefore the equations~(\ref{Eqn:WPWS:Majorant}) and (\ref{Eqn:GFMW:RatioFEAxeReel}) imply
\begin{equation}\label{Eqn:GFMW:MajorantSigma}
 M(x)\simeq\ex{h(x)}\ex{\omega_{\psi'/\pi}(x)}\sqrt{\sigma(x)},
\end{equation}
 where $M$ is the majorant of $\Hil(E)$.

Our aim is to express $\sqrt{\sigma(x)}$ in the form $\simeq\ex{\omega_m(x)}$ for a suitable, not necessarily positive $m$, in other words,
 to express $\dem\log\sigma(x)$ as a potential.
We reach the following general question: given an $f(x)$, is it possible to find a bounded $m$ such that $\omega_m(x)=f(x)$ for all $x\in\reel$?
Intuitively, if $f$ is Poisson integrable, we may force the stronger relation $\omega_m(z)=P_f(z)$ on the upper half-plane.
Then, $\dpar{y}P_f=\pi P_m=-\dpar{x}\cjh{P_f}$, where the tilde denotes the harmonic conjugate.
A look at the limit when $y\downarrow 0$ suggests to consider $m=-\cjh{f'}/\pi$, where $\cjh{f'}$ is the Hilbert transform of $f'$.

In our case $f$ shall be a smoothing of $\dem\log\sigma(x)$, whose exact definition shall be given later.
By the axiom (\ref{Eqn:MCA:CroissanceSommets1}), $f(x)-f(x')$ is expected to be $O(|x-x'|^{1-\eps})$ uniformly in $x$ when $|x-x'|\to\infty$. 
In such circumstances the magnitude of $\cjh{f'}$ is controlled by the ones of $f'$ and $f''$, as shown by the following lemma. 

\begin{lem}
 Consider a twice differentiable function $f$ on $\reel$ satisfying $\nm{f'}<\infty$, $\nm{f''}_\infty<\infty$, and 
  $|f(x)-f(x')|\leq A|x-x'|^{1-\eps}$ uniformly in $x$ when $|x-x'|\geq R$ for a certain $A>0$ and a certain $R>1$.
 Then $\cjh{f'}$ exists and satisfies
  $$\pi\nm{\cjh{f'}}_\infty\leq 2\nm{f''}_\infty+2\log R\ \nm{f'}_\infty+\frac{2A}{R^\eps}\left(1+\Inv{\eps}\right).$$
\end{lem}

\begin{pr}
 Recall that
  \begin{equation}\label{Eqn:GFMW:TransformeeHilbert}
   \pi\cjh{f'}(x)=\vp_{-\infty}^\infty\frac{f'(t)}{x-t}\dif{t}=\int_0^\infty\frac{f'(x-t)-f'(x+t)}{t}\dif{t}
  \end{equation}
 if this last improper integral is well-defined.
 The mean value theorem implies
  $$\int_0^1\left|\frac{f'(x-t)-f'(x+t)}{t}\right|\dif{t}\leq 2\nm{f''}_\infty.$$
 Moreover,
  $$\int_1^R\left|\frac{f'(x-t)-f'(x+t)}{t}\right|\dif{t}\leq 2\nm{f'}_\infty\int_1^R\Inv{t}\dif{t}=2\log R\ \nm{f'}_\infty.$$
 Finally, an integration by parts yields
  $$\int_R^\infty\frac{f'(x-t)}{t}\dif{t}=\int_R^\infty\frac{f(x)-f(x-t)}{t^2}\dif{t}-\frac{f(x)-f(x-R)}{R}.$$
 By the assumption this last quantity is bounded in absolute value by $\Kfrac{A}{R^\eps}\left(1+\Inv{\eps}\right)$.
 A similar conclusion holds for $\Kint_R^\infty\frac{f'(x+t)}{t}\dif{t}$.
 Hence, $\cjh{f'}$ exists and satisfies the desired estimate.
\end{pr}


Our previous discussion may be formalized as follows.

\begin{lem}
 Under the assumptions of the previous lemma, there exist $a,c\in\reel$ such that
  $f(x)=\omega_{-\cjh{f'}/\pi}(x)+ax+c$
 for $x\in\reel$.
 More generally, for $z$ in the upper half-plane
  $P_f(z)=\omega_{-\cjh{f'}/\pi}(z)+a\Re z+c.$
\end{lem}

\begin{pr}
 Notice that $\omega_{\cjh{f'}/\pi}(z)$ is well-defined, since $\cjh{f'}$ is bounded.
 It satisfies $\dpar{y}\omega_{\cjh{f'}/\pi}(z)=P_{\cjh{f'}}(z)$ for all $z\in\cpl^+$.
 The Poisson transform of $f'$ is also well-defined, since $f'$ is bounded.
 Moreover an integration by parts shows that the harmonic conjugate
  $$\cjh{P_{f'}}(z)=\Inv{\pi}\int_{-\infty}^\infty-\Re\left(\Inv{t-z}\right)f'(t)\dif{t}$$
 is well-defined and bounded, though it is not necessarily absolutely convergent.
 Let us show that $\cjh{P_{f'}}=P_{\cjh{f'}}$, even if the last integrand and the one in (\ref{Eqn:GFMW:TransformeeHilbert}) were not in $L^1$.
 In fact,
 $$\pi\cjh{P_{f'}}(z)=\!\!\int_{-\infty}^\infty\frac{x-t}{(x-t)^2+y^2}f'(t)\dif{t}=\!\!\int_0^\infty\frac{f'(x-t)-f'(x+t)}{t}\,\frac{t^2}{t^2+y^2}\dif{t}.$$
 Consequently, for $y>0$
  $$|\cjh{P_{f'}}(x+\eul y)-\cjh{f'}(x)|\leq \frac{2}{\pi}\,\nm{f''}_\infty\int_0^\infty\frac{y^2}{t^2+y^2}\dif{t}=\nm{f''}_\infty\cdot y,$$
 and hence $\lim_{y\downarrow 0}\cjh{P_{f'}}(x+\eul y)=\cjh{f'}(x)$ for all $x\in\reel$.
 Since $\cjh{P_{f'}}$ and $P_{\cjh{f'}}$ are bounded, the maximum principle \cite[lec.6 th.1]{Levin} implies $\cjh{P_{f'}}=P_{\cjh{f'}}$, as claimed.
 In addition, for $\Im z>0$
  $$\pi\cjh{P_{f'}}(z)=\int_{-\infty}^\infty-\Re\left(\Inv{(t-z)^2}\right)f(t)\dif{t}=\int_{-\infty}^\infty-\dpar{y}\Im\left(\Inv{t-z}\right)f(t)\dif{t}.$$
 The dominated convergence theorem%
\footnote{%
 Indeed, a dominator is easily derived from the mean value theorem, for instance, a constant times
  $$\frac{|t|^{1-\eps}+1}{(t-x)^2+(y/2)^2}.$$}
 and the well-definiteness of $P_f$ ensure that this last integral is equal to $-\dpar{y}P_f(z)$.
 In total, $\dpar{y}\omega_{\cjh{f'}/\pi}(z)=-\dpar{y}P_f(z)$ when $\Im z>0$.
 Since $\omega_{\cjh{f'}/\pi}$ and $P_f$ are harmonic, the result on $\cpl\setminus\reel$ follows.
 Since in addition $\omega_{\cjh{f'}/\pi}$ and $f$ are continuous, the result on $\reel$ also follows.
\end{pr}

We now construct a family $\sgl{f_L}$ of smoothings of $\dem\log\sigma$ such that: 
 each $f_L$ satisfies the hypotheses of the lemmas with a common $C$ and $\eps$;
 $\sqrt{\sigma(x)}\simeq\ex{f_L(x)}$ for any $L$; 
 and $\nm{\cjh{f_L'}}_\infty$ is as small as desired when $L$ is suitably chosen.
As we have seen, this last condition is realized when $\nm{f_L'}_\infty$ and $\nm{f_L''}_\infty$ are controlled.
For this reason, smoothings of the form $\rho*\dem\log\sigma$, where $\rho\in C^\infty(\reel)$ is compactly supported, are not suitable.
In fact, the variations of the slopes of $\rho*\dem\log\sigma$ are controlled by the vertical distances between successive steps in the graph of $\dem\log\sigma$.
Though these distances are bounded, they are not arbitrarily small.

The idea of Lyubarskii and Seip \cite{LyuSeip} is to construct the desired smoothings by replacing $\dem\log\sigma$ with a polygonal line
 whose vertices, which lie on the graph of $\dem\log\sigma$, are arbitrarily distant.
The slopes of the resulting segments are then controlled by the equation (\ref{Eqn:MCA:CroissanceSommets1}).
By smoothing this polygonal line, the resulting second derivative is then controlled by the differences of such consecutive slopes.

Formally, for an arbitrary $L\in\nat^*$, let $p_L(x)$ be the polygonal line joining the vertices $(jL,\dem\log\sigma(jL))$, where
 $j\in\rel$, and let $f_L(x)=(\rho*p_L)(x)$.
We claim that, if $L$ is large enough,

\begin{enumerate}
 \item \label{Enu:GFMW:BonneApproximation}
  $|\dem\log\sigma(x)-f_L(x)|\lesssim 1$;
 \item \label{Enu:GFMW:DecroissanceUniforme}
  $|f_L(x)-f_L(x')|\leq 4C|x-x'|^{1-\eps}$ if $|x-x'|\geq 3L$, uniformly in $x$ and $L$;
 \item \label{Enu:GFMW:ControleDeriveeSeconde}
  $\nm{f_L''}_\infty\leq 4\nm{\rho}_\infty CL^{-\eps}$;
 \item \label{Enu:GFMW:DeriveeBornee}
  $\nm{f_L'}_\infty\leq 2CL^{-\eps}$.
\end{enumerate}

In the following argument, $I$ denotes the support of $\rho\in C^\infty(\reel)$.
Though it is not essential, we assume that $|I|$ is less than the third of the minimal length of the mountain bases, that
 $\rho\geq 0$, and that $\int_I\rho(t)\dif{t}=1$.

Let us check the condition~\ref{Enu:GFMW:BonneApproximation}.
Let $x\in\reel$ and $k$ be the index of the first mountain whose base intercepts $[x,\infty)$.
If $x$ is at distance at least $2L$ from the plateaux, then by the equation (\ref{Eqn:MCA:CroissanceSommets1})
 \begin{eqnarray*}
  |f_L(x)-\dem\log\sigma(x)|&\leq&\sup_{t\in I}|p_L(x-t)-\dem\log\sigma(x)|\\
   &\leq&\dem \max_{j=-M}^{M}|\log\eta_{k+j}-\log\eta_k|\ \lesssim \ 1,
 \end{eqnarray*}
where $M$ exceeds by one the maximal number of mountains in an interval of length $L$. 
Otherwise, $x$ is close to a plateau, so a linking condition (similar to the argument in Proposition~\ref{Prop:MCA:CroissanceSommets}) ensures that $|\log\eta_{k+j}|\simeq 1$ for all $j=-M,\ldots,M$.
The result follows.

Regarding the condition~\ref{Enu:GFMW:DecroissanceUniforme}, since $|x-x'|\geq 3L$, there exists a $jL$ between $x$ and $x'$ at distance at least 
 $L$ from them.
If $L$ is large enough, the equation (\ref{Eqn:MCA:CroissanceSommets1}) then yields
 $$|\log\sigma(\xi)-\log\sigma(jL)|\leq 2C|\xi-jL|^{1-\eps}$$
for $\xi=x$, $\xi=x'$, and any $\xi$ further from $jL$.
The convexity of the set $\ens{(\xi,y)}{|y-\log\sigma(jL)|\leq 2C|\xi-jL|^{1-\eps}}$ on both sides of the line $\xi=jL$ ensures that
 $$|p_L(\xi)-\log\sigma(jL)|\leq 2C|\xi-jL|^{1-\eps}$$
for such $\xi$, and similarly for $f_L(x)$.
A fortiori, $|f_L(x)-\log\sigma(jL)|\leq 2C|x-jL|^{1-\eps}$ and $|f_L(x')-\log\sigma(jL)|\leq 2C|x'-jL|^{1-\eps}$.
Since $jL$ is between $x$ and $x'$, it follows that $|f_L(x)-f_L(x')|\leq 4C|x-x'|^{1-\eps}$. 

Regarding \ref{Enu:GFMW:ControleDeriveeSeconde}, on the one hand if $x$ is far from the vertices of $p_L$, then clearly $f_L''(x)=0$.
On the other hand, if $x$ is close to a vertex joining, say, segments of slopes $m_0$ and $m_1$, an integration by parts and the
 equation~(\ref{Eqn:MCA:CroissanceSommets1}) show that
 $$|f_L''(x)|=|(\rho''*p_L)(x)|\leq\nm{\rho}_\infty\,|m_0-m_1|\leq 4\nm{\rho}_\infty CL^{-\eps}.$$

Finally, the condition~\ref{Enu:GFMW:DeriveeBornee} is proven in a similar way.
Obviously $f_L'(x)=p_L'(x)$ if $x$ is far from the vertices.
Otherwise, $x$ is close to a vertex and then $f_L'(x)$ is a convex combination of the slopes surrounding this last.
Therefore, $\nm{f_L'}_\infty=\nm{p_L'}_\infty\leq 2CL^{-\eps}$ by the equation (\ref{Eqn:MCA:CroissanceSommets1}).

We are ready to conclude.
The conditions~\ref{Enu:GFMW:DecroissanceUniforme}, \ref{Enu:GFMW:ControleDeriveeSeconde}, and~\ref{Enu:GFMW:DeriveeBornee} ensure that the first
 lemma applies with $A=4C$ and $R=3L$, in particular, that
  $$\nm{\cjh{f'}}_\infty\lesssim\frac{\log L+1}{L^\eps}.$$
By the last lemma, the condition~\ref{Enu:GFMW:BonneApproximation} becomes $|\dem\log\sigma(x)+\omega_{\cjh{f_L'}/\pi}(x)-ax-c|\lesssim 1$.
Substituting the result in the relation~(\ref{Eqn:GFMW:MajorantSigma}),
 \begin{equation}\label{Eqn:GFMW:AdditionPotentiel}
  M(x)\simeq\exp\left(h(x)+ax+c+\omega_{(\psi'-\cjh{f_L'})/\pi}(x)\right).
 \end{equation}
We may choose $L$ so large that $\nm{\cjh{f_L'}}_\infty<\nm{\psi'}_\infty$.
Then $\psi'-\cjh{f_L'}$ is positive and comparable with $1$.
Therefore, the corollary of Proposition~\ref{Prop:GFMW:CasFacile} persists in the difficult case where $\ph'$ is not bounded.
In total, we have proven:

\begin{thm} \label{Thm:GFMW:GFMW}
 Let $M$ be the majorant of a given MC-space.
 There exists a measurable, positive function $m\simeq 1$ and a real-entire function $g$ such that
  $$M(x)\simeq\ex{g(x)}\ex{\omega_m(x)}$$
 when $x$ varies in $\reel$.
\end{thm}

\section{\!\!\!\!Generic form of the weighted PW-spaces}

Following the approach of Lyubarskii and Seip \cite{LyuSeip}, we shall compute the explicit form of the weighted PW-spaces.
Our main computational tool is the following proposition, based on this simple observation: \emph{in a mountain chain the integrals of the mountains are $\simeq 1$}.
This last assertion follows easily from the properties of the mountains, namely: they have Poissonian shapes, they consist of two sides and a summit bounded away from the extremities (in the horizontal direction), their summits have minimal height, and their bases have minimal length.

\begin{prop}
 Let $\gamma(x)$ be a mountain chain.
 Suppose that $f(z)$ is analytic in the strip $|\Im z|<\eps$, where $\eps>0$.
 Then, for $p>1$
  $$\int_{-\infty}^\infty|f(t)|^p\gamma(t)\dif{t}\lesssim\sup_{|y|<\eps}\int_{-\infty}^\infty|f(t+\eul y)|^p\dif{t}.$$
\end{prop}

\begin{pr}
 Since $p>1$, the mean value property and Jensen's inequality \cite[th.3.3]{Rudin} imply
  $$|f(t)|^p\leq\Inv{2\pi}\int_0^{2\pi}|f(t+r\exi{\theta})|^p\dif{\theta}$$
 for $0<r<\eps$.
 Consequently,
  $$|f(t)|^p\leq\Inv{\pi\eps^2}\int_{|z-t|<\eps}|f(z)|^p\dif{A}(z),$$
 where $\dif{A}(z)$ is the element of area.
 In particular, if $[a,b]$ is the base of a mountain,
  \begin{eqnarray*}\int_a^b|f(t)|^p\gamma(t)\dif{t}
   &\leq&   \Inv{\pi\eps^2}\int_{-\eps}^\eps\int_{a-\eps}^{b+\eps}|f(x+\eul y)|^p\dif{x}\dif{y}\int_a^b\gamma(t)\dif{t}\\
   &\simeq& \int_{-\eps}^\eps\int_{a-\eps}^{b+\eps}|f(x+\eul y)|^p\dif{x}\dif{y},
  \end{eqnarray*}
 independent of the mountain.
 Since in addition the mountain bases have length $\simeq 1$, while $\gamma(t)\simeq 1$ when $t$ is under a plateau, we conclude
  $$\int_{-\infty}^\infty|f(t)|^p\gamma(t)\dif{t}\lesssim \int_{-\eps}^\eps\int_{-\infty}^{\infty}|f(x+\eul y)|^p\dif{x}\dif{y}
  \lesssim\sup_{|y|<\eps}\int_{-\infty}^\infty|f(x+\eul y)|^p\dif{x}.$$
\end{pr}

On the one hand, given an MC-space $\Hil$ of majorant $M$, by Theorem~\ref{Thm:GFMW:GFMW} the restriction of $M$ to the real axis admits a representation
 $M\simeq \ex{g}\ex{\omega_m}$ for a real-entire $g$ and a measurable $m\simeq 1$.
On the other hand, we have shown in Section~\ref{Sec:APE} that $\ex{g}PW(m)$ is a weighted PW-space whose majorant-weight
 is $\simeq\ex{g}\ex{\omega_m(x)}$.
This suggests to compare $\Hil$ with $\ex{g}PW(m)$.

The following lemma could be stated as a corollary of the above proposition.
Notice that, if the concerned majorant admits several representation $\ex{g}\ex{\omega_m}$, any of them may be used.
Similarly, any choice of $E$ may be used.

\begin{lem}
 Let $\Hil=\Hil(E)$ be an MC-space, where $E$ is a considered HB-function satisfying the mountain chain axioms.
 Assume that the majorant of $\Hil$ is comparable with $\ex{g}\ex{\omega_m}$ on the real axis, where $g$ is real-entire and $m\simeq 1$ is measurable.
 Then, for any $f\in\ex{g}PW(m)$
  $$\int_{-\infty}^\infty\frac{|f(t)|^2}{|E(t)|^2}\dif{t}\lesssim\nm{f}_{\ex{g}PW(m)}^2.$$
\end{lem}

\begin{pr}
 Let $\ph$ be the phase of $E$.
 Suppose $f\in\ex{g}PW(m)$.
 The equation~(\ref{Eqn:WPWS:Majorant}) and the hypothesis on the majorant imply
  $$\int_{-\infty}^\infty\frac{|f(t)|^2}{|E(t)|^2}\dif{t}\simeq\int_{-\infty}^\infty|f(t)|^2\ex{-2g(t)}\ex{-2\omega_m(t)}\ph'(t)\dif{t}.$$
 In view of applying the last proposition (with $p=2$ and $\eps=1$), let us replace $\omega_m(t)$ with $\omega_m(t+\eul)$ in this last relation.
 This replacement does not break the equivalence, because
  $|\omega_m(t+\eul)-\omega_m(t)|\leq\pi\nm{m}_\infty.$
 Doing so,
  \begin{eqnarray*}\int_{-\infty}^\infty\frac{|f(t)|^2}{|E(t)|^2}\dif{t}
   &\simeq& \int_{-\infty}^\infty|f(t)|^2\ex{-2g(t)}\ex{-2\omega_m(t+\eul)}\ph'(t)\dif{t}\\
   &\lesssim& \sup_{|y|< 1}\int_{-\infty}^\infty|f(t+\eul y)|^2\ex{-2\Re g(t+\eul y)}\ex{-2\omega_m(t+\eul y+\eul)}\dif{t}\\
   &\simeq& \sup_{|y|< 1}\int_{-\infty}^\infty|f(t+\eul y)|^2\ex{-2\Re g(t+\eul y)}\ex{-2\omega_m(t+\eul y)}\dif{t}.
  \end{eqnarray*}
 Observe that $f(z)\ex{-g(z)}{\ex{-\omega_m(z)-\eul\cjh{\omega}_m(z)}}$ is analytic in the upper half-plane, continuous
  up to the real axis and, by the corollary of Proposition~\ref{Prop:APE:Lifting}, has exponential type $0$.
 The Plancherel--P\'olya theorem \cite[lec.7 th.4]{Levin} then implies that the right hand side is equal to $\nm{f}_{\ex{g}PW(m)}^2$, as claimed.
\end{pr}

A natural question arises: is it true that $\ex{g}PW(m)\subseteq\Hil(E)$ under the hypotheses of the lemma?
In general it is not clear, but for the \emph{specific} $m$ and $g$ constructed in the proof of Theorem~\ref{Thm:GFMW:GFMW} the answer is
 positive.

In the trivial case where the zeroes of $E$ are bounded away from the real axis, $m$ and $g$ are specified by the
 proposition~\ref{Prop:GFMW:CasFacile}, so $m=\ph'/\pi$.
This last proposition and the definition of de Branges space immediately yield $\Hil(E)=\ex{g}PW(m)$, as desired.

In the other case, $m=(\psi'-\cjh{f_L'})/\pi$ and $g(z)=h(z)+az+c$ are specified by the equation~(\ref{Eqn:GFMW:AdditionPotentiel}).
The second lemma of Theorem~\ref{Thm:GFMW:GFMW} asserts that $\omega_{-\cjh{f_L'}}(z)+a\Re z+c=P_{f_L}(z)$.
Therefore, the condition 
 \begin{equation}\label{Eqn:GFWPWS:NevanlinnaBasique}
  f(z)\ex{-g(z)}\ex{-\omega_m(z)-\eul\cjh{\omega}_m(z)}\in\Nev_0^\pm
 \end{equation}
 is equivalent to $f(z)\ex{-h(z)}\ex{-\omega_{\psi'/\pi}(z)-\eul\cjh{\omega}_{\psi'/\pi}(z)}\in\Nev_0^\pm$.
Recall that $$\ex{\Re h(z)}\ex{\omega_{\psi'/\pi}(z)}\simeq |F(z)|$$ in the upper half-plane, where $F$ is obtained from $E$ by 
 lowering some of its zeroes.
Thus, the last condition is equivalent to
 $f(z)/F(z),\,f^*(z)/F(z)\in\Nev_0^+,$
and hence to $f(z)/E(z),\,f^*(z)/E(z)\in\Nev_0^+$, as shown by the proposition~\ref{Prop:GFMW:InterchangementEF}.
This last condition and the lemma imply that 
 $\ex{g}PW(m)\subseteq\Hil(E)$.

What else if $\Hil(E)$ is indeed a weighted PW-space?
By the previous chain of equivalences, if $f\in\Hil(E)$, then the relation~(\ref{Eqn:GFWPWS:NevanlinnaBasique}) is satisfied.
By the choice of $g$ and $m$, $\Hil(E)$ and $\ex{g}PW(m)$ have comparable majorants, and hence their respective norms are equivalent.
Therefore $f\in\ex{g}PW(m)$, that is, $\ex{g}PW(m)=\Hil(E)$.
In summary,

\begin{thm}\label{Thm:GFWPWS:GFWPWS}
 If $\Hil$ is an MC-space, then there exists a real-entire $g(z)$ and a measurable, positive $m(x)\simeq 1$ such that
  $\ex{g}PW(m)\subseteq\Hil$, while these two spaces have majorant $\simeq \ex{g(x)}\ex{\omega_m(x)}$ on the real axis.
 In particular, if $\Hil$ is a weighted Paley--Wiener space, $\ex{g}PW(m)=\Hil$ (equality with equivalence of norms).
\end{thm}

Notice that an inclusion of the form $\ex{g}PW(m)\subseteq\Hil$, where $\Hil$ is a weighted PW-space, does not  necessarily imply $\ex{g}PW(m)=\Hil$, since
 the involved $g$ and $m$ are not necessarily those constructed from $\Hil$ in the previous work.
\TroisEtoiles
From an $E$ satisfying the mountain chain axioms, we have constructed a real-entire $g$ and a measurable $m\simeq 1$
 such that $\ex{g}PW(m)\subseteq\Hil(E)$.
Let us clarify the relationship between $E$ and $\ex{g}\ex{\omega_m}$.

By lowering some zeroes of $E$ we have defined an HB-function $F$ satisfying
 $$|F(x+\eul|y|)|=\ex{\Re h(z)}\ex{\omega_{\psi'/\pi}(z)},$$
where $h$ is real-entire and $\psi'(x)\simeq 1$.
The contribution of the zeroes of $E$ in the critical strip has been retained by a step function $\log\sigma(x)$, replaced with a smoothed
 broken line $f_L$ such that $|\log\sigma(x)-f_L(x)|\lesssim 1.$

This last inequation implies
 $|P_{\log\sigma}(z)-P_{f_L}(z)|\lesssim 1$
for $\Im z>0$.
In addition, the second lemma of Theorem~\ref{Thm:GFMW:GFMW} implies
 $P_{f_L}(z)=\omega_{-\cjh{f_L'}/\pi}(z)+a\Re z+c,$
so
 $$\ex{P_{\log\sigma}(z)}\simeq\ex{\omega_{-\cjh{f_L'}/\pi}(z)}\ex{a\Re z+c}.$$
Therefore, the definitions of $g$ and $m$ given in the equation~(\ref{Eqn:GFMW:AdditionPotentiel}) yield
 $$|F(z)|\ex{P_{\log\sigma}(z)}\simeq\ex{\Re g(z)}\ex{\omega_m(z)}$$
for $\Im z>0$.
In other words, by the equation (\ref{Eqn:GFMW:RatioHBFunctions}), 
 \begin{equation}\label{Eqn:GFWPWS:LienEgm}
  |E(z)|\ex{P_{\log\sigma}(z)}\simeq\min(1,|z-\lambda_0(z)|)\ex{\Re g(z)}\ex{\omega_m(z)},
 \end{equation}
where $\lambda_0(z)$ is the zero of $E$ in the critical strip closest to $z$.

Notice that this last equation relates the majorant of $\Hil(E)$ with the majorant of $\ex{g}PW(m)$ for $\Im z\gg 0$.
Indeed, denoting the former by $M_{\Hil}$ and the latter by $M_{\pw}$, the equation~(\ref{Eqn:WPWS:NoyauReproduisant}) yields
 $$M_{\Hil}(z)^2=\frac{|E(z)|^2}{4\pi\Im z}(1-|E^*(z)/E(z)|^2),$$
while the relation~(\ref{Eqn:APE:MajorantDemiPlan}) yields
 $M_{\pw}(z)^2\simeq\kfrac{\ex{2g(z)}\ex{2\omega_m(z)}}{\Im z}.$
Therefore, (\ref{Eqn:GFWPWS:LienEgm}) is equivalent to
 $$M_{\Hil}(z)\ \ex{P_{\log\sigma}(z)}\simeq M_{\pw}(z)\ \sqrt{1-|E^*(z)/E(z)|^2}$$
for $\Im z\gg 0$.
In the case where $\Hil(E)=\ex{g}PW(m)$, we obtain
 $$\ex{P_{\log\sigma}(z)}\simeq \sqrt{1-|E^*(z)/E(z)|^2}$$
for $\Im z\gg 0$.

As observed by Lyubarskii and Seip \cite[p.1004 rem.2]{LyuSeip}, in the special circumstances where $\sigma(x)\to 0$ as $x\to+\infty$ (resp.\@ $x\to-\infty$), the above
 relation forces the zeroes of $E$ in a right half-plane (resp.\@ left half-plane) to lie in a horizontal strip.
Let us treat the case where $\lim_{x\to +\infty}\sigma(x)=0$, the other case being similar.
Suppose by contradiction that there exists a sequence $\sgl{\xi_k-\eul\eta_k}$ of zeroes of $E$ such that $\eta_k\to\infty$, while $\xi_k$ is bounded below.
Then, $\ex{P_{\log\sigma}(\xi_k+\eul\eta_k)}\simeq 1$, that is, $|P_{\log\sigma}(\xi_k+\eul\eta_k)|\lesssim 1$.
In other words,
 $$\int_{-\infty}^\infty\frac{-\log\sigma(\eta_k t+\xi_k)}{t^2+1}\dif{t}\lesssim 1.$$
This contradicts Fatou's lemma, since $\lim_{k\to\infty}-\log\sigma(\eta_k t+\xi_k)=\infty$ for all positive $t$.

Let us add that in these circumstances any zero of $E$ whose real part is sufficiently large lies in the critical strip, and hence corresponds to
 the summit of a mountain.
Otherwise, there is a sequence of zeroes of $E$ whose $x$-coordinates go to $+\infty$, while their $y$-coordinates are $\simeq 1$
 in absolute value.
Eventually, each such zero would lie under a mountain and have a contribution $\simeq 1$ to $\ph'$ at the edge of this last.
This contradicts the axiom~\ref{Axi:MCA:ChaineMontagnes}.

Unfortunately the previous technique fails for the more general situation where $\liminf_{x\to+\infty}\sigma(x)=0$.

\section{\!\!\!\!Characterization of weighted PW-spaces}\label{Sec:NSC}

Which supplemental conditions an MC-space must satisfy in order to be a weighted PW-space?
As we shall see, a convenient condition may be obtained from the study of complete interpolating sequences in weighted PW-spaces.
Such sequences in classical PW-spaces have been characterized geometrically by Pavlov \emph{et al.} \cite{Pavlov}.
Following the line of Lyubarskii and Seip, one may lift the result to the weighted case.

In the sequel, given a countable $\Lambda\subset\cpl$ and a measurable $m\simeq1$, the discrete version of $PW(m)$ on $\Lambda$ is denoted by $\pw(\Lambda,m)$.
It consists of all sequences $\fml{c_\lambda}{\lambda\in\Lambda}{}$ whose norm
 $$\nm{\sgl{c_\lambda}}_{\pw(\Lambda,m)}=\left(\sum_{\lambda\in\Lambda}|c_\lambda|^2\ex{-2\omega_m(\lambda)}(1+|\Im\lambda|)\right)^{1/2}$$
is finite.
For $\tau>0$, the classical spaces $PW(\tau)$ and $\pw(\Lambda,\tau)$ are denoted by $L^2_{\pi\tau}$ and $\ell^2_{\pi\tau}(\Lambda)$ respectively.
Moreover, $L^2_{\pi\tau}[\Sigma]$ denotes the annihilator in $L^2_{\pi\tau}$ of a given $\Sigma\subseteq\cpl$.
It is a subspace of $L^2_{\pi\tau}$, but obviously not a de Branges space in general.

With the above definitions a sequence $\Lambda$ is said to be \emph{interpolating} \cite{LyuSeip} for $PW(m)$ if for each $\sgl{c_\lambda}\in\pw(\Lambda,m)$ there exists
 an $f\in PW(m)$ such that $f(\lambda)=c_\lambda$ for all $\lambda\in\Lambda$.
It is \emph{complete interpolating} if in addition $f$ is unique.

It is then natural to propose similar definitions for $\ex{g}PW(m)$, where the eligible interpolation data are taken 
 in $\ex{g}\pw(\Lambda,m)$. 
One then sees that a sequence $\Lambda$ is [complete] interpolating for $\ex{g}PW(m)$ if and only if it is so 
 for $PW(m)$.

We now complete the statement of Proposition~\ref{Prop:APE:Lifting} with its discrete analogue, whose proof is 
 straightforward.
 
\begin{thm}\label{Thm:NSC:Lifting}
 Given a real-entire $g$ and a measurable $m\simeq1$, assume $\tau>\sup m$.
 Let $E_{\tau-m}$ be given by Proposition~\ref{Prop:APE:HBCanonique}, and let us denote its zero set by $\Sigma$.
 Then, $f\mapsto f\ex{-g}E_{\tau-m}$ is a bijection from $\ex{g}PW(m)$ to $L^2_{\pi\tau}[\Sigma]$ with equivalence of norms.
 Moreover, given a countable $\Lambda\subset\adh{\cpl^+}$, $\sgl{c_\lambda}\mapsto\sgl{\ex{-g(\lambda)}E_{\tau-m}(\lambda)c_\lambda}$ is a bijection from
  $\ex{g}\pw(\Lambda,m)$ to $\ell^2_{\pi\tau}(\Lambda)$ with equivalence of norms.
\end{thm}

Extending the definition of complete interpolating sequence in the obvious way, it is then immediate from the following commutative squares
 that \emph{a sequence $\Lambda$ is complete interpolating for $\ex{g}PW(m)$ if and only if it is complete interpolating for $L^2_{\pi\tau}[\Sigma]$.}
Here, the well-definiteness of the left (resp.\@ right) side of the square is ensured by the fact that $\Lambda$ is complete interpolating for $\ex{g}PW(m)$
 (resp.\@ $L^2_{\pi\tau}[\Sigma]$).

\begin{equation*}
 \begin{array}{ccccccc}
  \ex{g}PW(m)&\leftrightarrow&L^2_{\pi\tau}[\Sigma]&\hspace{1.5cm}&\ex{g}PW(m)&\leftrightarrow&L^2_{\pi\tau}[\Sigma]\\
  \uparrow&&\exists!\uparrow&&\exists!\uparrow&&\uparrow\\
  \ex{g}\pw(\Lambda,m)&\leftrightarrow&\ell^2_{\pi\tau}(\Lambda)&&\ex{g}\pw(\Lambda,m)&\leftrightarrow&\ell^2_{\pi\tau}(\Lambda)
 \end{array}
\end{equation*}\Rtr

A classical result of Beurling asserts that, \emph{if a real sequence $\Gamma$ is separated and satisfies $D^+(\Gamma)<\tau$, then it is
 interpolating for $L^2_{\pi\tau}$}, where
 $$D^+(\Gamma)=\lim_{r\to\infty}\sup_{a\in\reel}\Inv{r}\#(\Gamma\cap[a,a+r]).$$
In particular, since $D^+(\Sigma+\eul)\leq\nm{\tau-m}_\infty$ (as shown by the construction of $E_{\tau-m}$ in Proposition~\ref{Prop:APE:HBCanonique}),
 the sequence $\Sigma$ is interpolating for $L^2_{\pi\tau}$.
This yields:

\begin{cor}
 A sequence $\Lambda$ in the closed upper half-plane is complete interpolating for $\ex{g}PW(m)$ if and only if $\Lambda\cup\Sigma$ is complete
  interpolating for $L^2_{\pi\tau}$.
\end{cor}

\begin{pr}
 Obviously, if $\Lambda\cup\Sigma$ is complete interpolating for $L^2_{\pi\tau}$, then $\Lambda$ is complete interpolating for $L^2_{\pi\tau}[\Sigma]$,
  equivalently, for $\ex{g}PW(m)$.
 Conversely, suppose $\Lambda$ is complete interpolating for $L^2_{\pi\tau}[\Sigma]$.
 Let $\sgl{c_\alpha}\in\ell^2_{\pi\tau}(\Lambda\cup\Sigma)$.
 Since $\Sigma$ is interpolating, there exists an $f\in L^2_{\pi\tau}$ satisfying $f(\sigma)=c_\sigma$ for all $\sigma\in\Sigma$.
 Observe that $\Lambda$ is interpolating for $L^2_{\pi\tau}$.
 A classical theorem of Plancherel and P\'olya \cite[lec.20 th.3]{Levin} then implies that $\sgl{f(\lambda)}\in\ell^2_{\pi\tau}(\Lambda)$.
 In particular, since $\Lambda$ is interpolating, there exists a $g\in L^2_{\pi\tau}[\Sigma]$ satisfying $g(\lambda)=c_\lambda-f(\lambda)$ for
  all $\lambda\in\Lambda$.
 Then, $(f+g)(\alpha)=c_\alpha$ for all $\alpha\in\Lambda\cup\Sigma$, as desired.
 Finally, $f+g$ is the only solution of this last interpolation problem: its difference with any possible solution
  vanishes, since it is in $L^2_{\pi\tau}[\Sigma]$, while $\Lambda$ is complete interpolating for this last space.
\end{pr}
%
\kTroisEtoiles
We now make use of the following, classical theorem of de Branges \cite[th.22]{deBranges}.

\begin{prop}\label{Prop:CWPWS:deBranges}
 Let $E$ be an Hermite--Biehler function of phase $\ph$.
 Let $\Lambda_\alpha$ be the solution set of
  $\ph(x)\equiv\alpha\ (\mod\pi),$
 $0\leq\alpha<\pi$.
 If $\exi{\alpha}E(z)-\exni{\alpha}E^*(z)\notin\Hil(E)$, then the reproducing kernels $\sgl{k_\lambda}_{\lambda\in\Lambda_\alpha}$ constitute an orthogonal
  basis in $\Hil(E)$.
 Furthermore, the relation $\exi{\alpha}E(z)-\exni{\alpha}E^*(z)\in\Hil(E)$ may occur for at most one exceptional $\alpha$.
\end{prop}

Let us show \cite[p.1003]{LyuSeip} that, if $\Hil(E)$ is a weighted PW-space, the exceptional $\alpha$ does not exist.
This may be deduced the following, more general proposition.

\begin{prop}
 Let $E$ be a considered HB-function of phase $\ph$ satisfying the mountain chain axioms. 
 Let $\sigma(x)=\min(1,\eta(x))$,
  where $\xi(x)-\eul\eta(x)$ denotes the zero of $E$ closest to $x$.
 If $\sin^2\ph(x)$ is integrable, then so is $\sigma(x)$.
\end{prop}

\begin{pr}
 Suppose $\int_{-\infty}^\infty\sin^2\ph(x)\dif{x}<\infty$.
 Then, for any choice of $\eps>0$, $|\sin\ph(x)|>\eps$ may only occur on a set of finite measure.
 It forces $\sigma(x)$ to tend to $0$ when $|x|\to\infty$, so $\ph(x)$ resembles more and more a stair when $|x|$ is large enough: eventually each
  ``horizontal'' panel of this stair must be at level near to an integral multiple of $\pi$, so $\sin\ph(x)$ remains small along it.

 Let $\sgl{a_k}$ be the solution set of $\sin\ph(x)=1$.
 Necessarily for $|a_k|$ large enough the $\sgl{(a_k,\ph(a_k))}$ intercept all the ``vertical'' panels of the aforementioned stair, maybe several times.
 Let $b_k$ be the smallest solution of $\sin\ph(x)=1/2$ larger than $a_k$.
 Then, the mean value theorem implies the existence of $r\in(a_k,b_k)$ such that
  $$\frac{\sin\ph(a_k)-\sin\ph(b_k)}{b_k-a_k}=\Inv{2(b_k-a_k)}=|\cos\ph(r)|\ph'(r)\leq\ph'(r).$$
 This last quantity is $\lesssim\kInv{\sigma(a_k)}$, since the mountain over $a_k$ has height $\kInv{\sigma(a_k)}$.
 Therefore, there exists a $c\in(0,1)$ such that $b_k-a_k\geq c\sigma(a_k)$ for all $k$ large enough (in absolute value).
 Repeating the argument with $b_k$ to the left of $a_k$, we conclude that $\sin\ph(x)\simeq 1$ when $|x-a_k|<c\sigma(a_k)$.
 In particular,
  $$\int_{-\infty}^\infty\sin^2\ph(x)\dif{x}\geq\sum_{|k|\ \mathrm{large}}\int_{|x-a_k|<c\sigma(a_k)}\sin^2\ph(x)\dif{x}\simeq \sum_{|k|\ \mathrm{large}}\sigma(a_k).$$
 Since all vertical panels of the stair are intercepted at least once, $\Kint_{|x|>R}\!\!\!\!\!\!\!\!\!\!\!\sigma(x)\dif{x}$ is comparable with this last sum
  for $R$ large enough.
 The result follows.
\end{pr}

\begin{cor}
 Let $E$ be a considered HB-function and $0\leq\alpha<\pi$.
 If $\Hil(E)$ is a weighted Paley--Wiener space, then $\exi{\alpha}E(z)-\exni{\alpha}E^*(z)\notin\Hil(E)$.
\end{cor}

\begin{pr}
 It suffices to show that $\sin(\ph(x)-\alpha)\notin L^2(\reel)$, where $\ph$ is the phase of $E$.
 Assume the converse, so in particular $\lim_{|x|\to\infty}\sigma(x)=0$.
 Since $\exi{\alpha}E(z)$ satisfies the mountain chain axioms and shares its zeroes with $E$, by the proposition
  $\sigma(x)\in L^1(\reel)$.

 Let $\zeta$ be a zero of $E$.
 Then $E(z)/(z-\zeta)\in\Hil(E)$.
 Since $\Hil(E)$ is a weighted PW-space, $\nm{E(x)/(M(x)(x-\zeta))}_2<\infty$, where $M(x)$ is given by the equation~(\ref{Eqn:WPWS:Majorant}).
 Therefore,
  $$\int_{-\infty}^\infty\Inv{(x^2+1)\ph'(x)}\dif{x}<\infty.$$
 Observe that for $|x|$ large, $\sigma(x)$ is close to $0$, and hence $\ph'(x)\simeq\sigma(x)$ when $x$ is 
  bounded away from the middle of its corresponding mountain (independent of the mountain).
 Consequently, there exists a set $L$, which consists of an infinite union of disjoint intervals of radii $\gtrsim 1$,
  such that
   $$\int_L\Inv{(x^2+1)\sigma(x)}\dif{x}<\infty.$$
 However, Schwarz' inequality implies
  $$\left(\int_L\Inv{(x^2+1)\sigma(x)}\dif{x}\right)^{1/2}\left(\int_L\sigma(x)\dif{x}\right)^{1/2}\geq\infty,$$
 and hence $\Kint_L\sigma(x)\dif{x}=\infty$, a contradiction.
\end{pr}

As a consequence, if $\Hil(E)=\ex{g}PW(m)$ is a weighted PW-space, the conclusion of
 Proposition~\ref{Prop:CWPWS:deBranges} holds for \emph{all} $\alpha\in[0,\pi)$.
In particular, $\Lambda_\alpha$ is complete interpolating for $\ex{g}PW(m)$.
By the corollary of Theorem~\ref{Thm:NSC:Lifting}, $\Lambda_\alpha\cup\Sigma$ is then complete interpolating for $L^2_{\tau\pi}$.

The following, geometric characterization was obtained by Pavlov \emph{et al.} \cite{Pavlov} and revisited by Lyubarskii and Seip \cite{LyuSeipA2} in a previous work.

\begin{prop}\label{Prop:NSC:LyubarskiiSeipClassique}
 A sequence $\Lambda\subset\cpl$ is complete interpolating for $L^2_{\pi\tau}$ if and only if the following conditions are satisfied:
  \begin{enumerate}
   \item
    $\Lambda$ is \emph{separated} for the metric $\rho(z,w)=\Kfrac{|z-w|}{1+|z-\conj{w}|}$ on $\cpl$, that is,
     $$\inf_{\dbl{\lambda,\lambda'\in\Lambda}{\lambda\neq\lambda'}}\rho(\lambda,\lambda')\gtrsim 1.$$
   \item
    $\Lambda$ satisfies the \emph{bilateral Carleson condition}, that is,
     $$\sup_{x\in\reel}\sum_{\dbl{\lambda\in\Lambda}{|\lambda-x|\leq R}}|\Im\lambda|\lesssim R.$$
   \item \label{Enu:CWPWS:ExponentialType}
    The \emph{generating product}
     $$S_\Lambda(z)=\lim_{R\to\infty}\prod_{\stackrel{\scriptstyle\lambda\in\Lambda}{|\lambda|\leq R}}\mbox{}\!\!\!'\,(1-z/\lambda)$$
    converges uniformly on each compact to a function of exponential type at most $\pi\tau$.
    (Here the apostrophe indicates that a possible factor indexed by $\lambda=0$ must be replaced with $z$.)
   \item\label{Enu:NSC:Muckenhoupt}
    The function $v(x)=|S_\Lambda(x)|^2/\dist{x}{\Lambda}^2$ satisfies the \emph{Muckenhoupt ($A_2$) condition}, that is,
     \begin{equation*}
      \Inv{|I|}\int_I v(x)\dif{x}\cdot\Inv{|I|}\int_I\Inv{v(x)}\dif{x}\lesssim 1
     \end{equation*}
    when $I$ varies among the finite intervals.
  \end{enumerate}
\end{prop}\mbox{}

\noindent In particular, if $\Hil(E)=\ex{g}PW(m)$, then $S_{\Lambda_\alpha\cup\Sigma}(z)$ converges uniformly on each compact to a function of exponential type at most $\pi\tau$, 
 while the weight $v(x)=|S_{\Lambda_\alpha\cup\Sigma}(x)|^2/\dist{x}{\Lambda_\alpha\cup\Sigma}^2$ satisfies the Muckenhoupt ($A_2$) condition.

Notice that $\Lambda_\alpha$ is precisely the zero set of the entire function $\exi{\alpha}E(z)-\exni{\alpha}E^*(z)$.
(The zeroes of this last function are all real, since $E$ is in the Hermite--Biehler class.)
Hence, $S_{\Lambda_\alpha\cup\Sigma}$ and $\ex{-g}E_{\tau-m}(\exi{\alpha}E-\exni{\alpha}E^*)$ have the same zeroes.
Let us show that these two functions are indeed equal, up to a multiplicative constant.

Let $\zeta$ be a zero of $E^\sharp$, where $E^\sharp\in\sgl{E,E^*}$.
Since $E^\sharp(z)/(z-\zeta)\in\Hil(E)$, Theorem~\ref{Thm:NSC:Lifting} implies that $\ex{-g(z)}E_{\tau-m}(z)E^\sharp(z)/(z-\zeta)\in L_{\pi\tau}^2$, so $\ex{-g(z)}E_{\tau-m}(z)E^\sharp(z)$ has exponential type $\leq\pi\tau$, while $$\ex{-g(x)}E_{\tau-m}(x)E^\sharp(x)/\sqrt{x^2+1}\in L^2(\reel).$$
A fortiori,
 $$\int_{-\infty}^\infty\frac{\log^+|\ex{-g(x)}E_{\tau-m}(x)E^\sharp(x)|}{x^2+1}\dif{x}<\infty,$$
so $\ex{-g}E_{\tau-m}(\exi{\alpha}E-\exni{\alpha}E^*)$ is in the Levinson--Cartwright class (a.k.a.\@ class~$C$, \cite[lec.16]{Levin}).
In addition, the penultimate statement of Proposition~\ref{Prop:APE:HBCanonique} implies 
 $$\limsup_{r\to\infty}\frac{\log|E_{\tau-m}(\eul r)|}{r}=\limsup_{r\to\infty}\frac{\log|E_{\tau-m}(-\eul r)|}{r},$$
and hence the analogue property for $\ex{-g}E_{\tau-m}(\exi{\alpha}E-\exni{\alpha}E^*)$ (instead of $E_{\tau-m}$) is straightforward,
 using the fact that $g$ is real-entire.
Finally, $\ex{-g}E_{\tau-m}(\exi{\alpha}E-\exni{\alpha}E^*)$ has at most a simple zero at $0$. 
The desired factorization follows \cite[lec.17 rem.2]{Levin}.

Thus, the square of 
 $$\kfrac{|\ex{-g(x)}E_{\tau-m}(x)(\exi{\alpha}E(x)-\exni{\alpha}E^*(x))|}{\dist{x}{\Lambda_\alpha\cup\Sigma}}$$
satisfies the Muckenhoupt($A_2$) condition. 
Notice that Proposition~\ref{Prop:APE:HBCanonique} and the equation~\ref{Eqn:WPWS:Majorant} imply
 $$|\ex{-g(x)}E_{\tau-m}(x)|\simeq\ex{-g(x)}\ex{-\omega_m(x)}\simeq\Inv{\sqrt{\ph'(x)}|E(x)|}$$
on the real axis.
In addition, $|\exni{\alpha}E^*(x)-\exi{\alpha}E(x)|=2|E(x)|\,|\sin(\ph(x)-\alpha)|.$
Finally, since $\ph'$ is comparable with a mountain chain, the distance between two successive elements in $\Lambda_\alpha$ is bounded above.
Hence,
  $$\dist{x}{\Lambda_\alpha\cup\Sigma}\simeq\dist{x}{\Lambda_\alpha}.$$
In conclusion, $\sin^2(\ph(x)-\alpha)/\ph'(x)\dist{x}{\Lambda_\alpha}^2$ must satisfy the Muckenhoupt ($A_2$) condition.
Using Proposition~\ref{Prop:NSC:LyubarskiiSeipClassique} integrally, we deduce:

\begin{prop}\label{Prop:CWPWS:ConditionZeros}
 Let $E$ be a considered Hermite--Biehler function of phase $\ph$ and let $\Lambda_\alpha\subset\reel$ be the zero set of $\sin(\ph(x)-\alpha)$.
 Suppose $\Hil(E)$ is a weighted Paley--Wiener space.
 Then, for all $\alpha\in[0,\pi)$, $\Lambda_\alpha$ is separated (for the Euclidean metric), while
  $$v(x)=\frac{\sin^2(\ph(x)-\alpha)}{\ph'(x)\dist{x}{\Lambda_\alpha}^2}$$
 satisfies the Muckenhoupt ($A_2$) condition.
\end{prop}

This last necessary condition creates a severe condition on the location of the zeroes of $E$.
For instance, denoting by $\xi_k-\eul\eta_k$ the zero corresponding to the $k$-th mountain, if the inferior limit of $\eta_k$ is going
 to zero too fast, it may contradict both the Mukenhoupt and the separation conditions.
The latter, because parts of the graph of $\ph$ would resemble stairs with vertical panels longer than $\pi$ and whose slopes are arbitrarily close
 to infinity.
Since these panels would intercept at least two points $(\ph^{-1}(\alpha+k\pi),\alpha+k\pi)$, $\Lambda_\alpha$ would not be separated.

Joined to a compatibility condition between $E$ and $\ex{g}\ex{\omega_m}$, the conclusion of the last proposition also gives a sufficient condition 
 \cite[th.4]{LyuSeip} for an MC-space to be a weighted PW-space (more precisely, to be the one constructed in Theorem~\ref{Thm:GFWPWS:GFWPWS}).

\begin{thm}
 Let $E$ be an HB-function of phase $\ph$ satisfying the mountain chain axioms, and $\Lambda_\alpha$ be the zero set
  of $\sin(\ph(x)-\alpha)$.
 Let $g$ and $m$ be given by Theorem~\ref{Thm:GFWPWS:GFWPWS}.
 Then, $\Hil(E)=\ex{g}PW(m)$ if and only if the following conditions are satisfied:
 Firstly, $\ex{-g(z)}\ex{-\omega_m(z)}E(z)$ is a (nonentire) function of exponential type on $\cpl$, which satisfies
    $$\int_{-\infty}^\infty\frac{\log^+|\ex{-g(x)}\ex{-\omega_m(x)}E(x)|}{x^2+1}\dif{x}<\infty.$$
 Secondly, for two $\alpha\in[0,\pi)$, $\Lambda_\alpha$ is separated, while
    $$v(x)=\frac{\sin^2(\ph(x)-\alpha)}{\ph'(x)\dist{x}{\Lambda_\alpha}^2}$$
 satisfies the Muckenhoupt ($A_2$) condition.
\end{thm}

\begin{pr}
 Necessity of the latter condition comes from Proposition~\ref{Prop:CWPWS:ConditionZeros}, while necessity of 
  the former is proven before, after Proposition~\ref{Prop:NSC:LyubarskiiSeipClassique}, using Theorem~\ref{Thm:NSC:Lifting}.
  
 For $\tau>\sup m$ large enough, let $E_{\tau-m}$ be given by Proposition~\ref{Prop:APE:HBCanonique} and let $\Sigma$ be
  its zero set.
 The first assumption implies that $\ex{-g}E_{\tau-m}(\exi{\alpha}E-\exni{\alpha}E^*)$ is an  entire function of
  exponential type $\leq\pi\tau$ on $\overline{\cpl^+}$, hence on $\cpl$, which satisfies the Levinson--Cartwright condition.
 Therefore \cite[lec.17 rem.2]{Levin}, it is equal to a constant times the generating product $S_{\Lambda_\alpha\cup\Sigma}$.
 In particular, this last satisfies the condition~\ref{Enu:CWPWS:ExponentialType} of
  Proposition~\ref{Prop:NSC:LyubarskiiSeipClassique}.
 In addition, $\Lambda_\alpha\cup\Sigma$ is clearly separated.
 Finally, $v(x)\simeq |S_{\Lambda_\alpha\cup\Sigma}(x)|^2/\dist{x}{\Lambda_\alpha}^2$ by the argument preceding
  Proposition~\ref{Prop:CWPWS:ConditionZeros}.
 Therefore, Proposition~\ref{Prop:NSC:LyubarskiiSeipClassique} implies that $\Lambda_\alpha\cup\Sigma$ is complete
  interpolating in $L^2_{\tau\pi}$, and hence, $\Lambda_\alpha$ is complete interpolating in $\ex{g}PW(m)$, as
  shown by the corollary of Theorem~\ref{Thm:NSC:Lifting}.
  
 This last conclusion holds for the two given $\alpha$.
 For at least one of them, by Proposition~\ref{Prop:CWPWS:deBranges}, $\Lambda_\alpha$ is complete interpolating
  also in $\Hil(E)$.
 Since such a sequence is in particular sampling, it follows that the norm of $f$ in $\Hil(E)$ is comparable with $\nm{f\ex{-g}\ex{-\omega_m}}_2$.  
 Moreover, by the proof of Theorem~\ref{Thm:GFWPWS:GFWPWS}, the condition $f/E, f^*/E\in\Nev_0^+$ is equivalent to
  $f\ex{-g}\ex{-\omega_m-\eul\cjh{\omega}_m}\in\Nev_0^\pm$.
 Therefore, $\Hil(E)=\ex{g}PW(m)$, as desired.
\end{pr}

\end{document}